\documentclass[a4paper,12pt,titlepage]{article}
\usepackage[latin1]{inputenc}
\usepackage[french]{babel}
\usepackage[dvips]{graphics}
\usepackage{amssymb,amsfonts}
\usepackage[T1]{fontenc}
\usepackage{amsmath}
\usepackage{amsthm}
\newcounter{def}
\newcounter{rem}

\newtheorem{lem}{Lemme}
\newtheorem{theo}{Th\'eor\`eme}

\newtheorem{cor}{Corollaire}
\newtheorem{prop}{Proposition}
\newtheorem{conj}{Conjecture}
\newcommand{\Ima}{\mathrm{\,Im\,}}
\newcommand{\Reel}{\mathrm{\,Re\,}}
\newcommand{\mc}{\mathcal }
\newcommand{\mf}{\mathfrak }
\newcommand{\mb}{\mathbb }

\def\build#1_#2^#3{\mathrel{\mathop{\kern 0pt#1}\limits_{#2}^{#3}}}
\newcommand{\ilbr}{\lbrack\kern-1.6pt\lbrack}
\newcommand{\irbr}{\rbrack\kern-1.6pt\rbrack}
\newcommand{\prim}{\begin{array}{l}\!\!\!\!\! \prime\\ \ \end{array}}

\newcommand{\bg}{\bigskip}

\newcommand{\Cl}{\mathrm{ Cl }}

\newcommand{\vol}{\mathrm{ Vol }}
\newcommand{\Tra}{\mathrm{Tr}}
\newcommand{\tra}{\mathrm{tr}}

\newcommand{\dis}{\displaystyle}

\newcommand{\sign}{\mathrm{\,sign}}
\begin{document}

\begin{center} \LARGE{Sommes de Dedekind associées à un corps de nombres totalement réel. }\\ \bg

\Large{Pierre Charollois}\end{center}\bg\bg\bg

\tableofcontents

\section{Introduction}

Soit $F$ un corps de nombres totalement réel de nombres de classes
1. On note $n=[F:\mb Q]$ son degré, $\mc O_F$ l'anneau des entiers
de $F,$ $\mc H$ le demi-plan de Poincaré, et on fixe un plongement
réel $\iota_j$ de $F.$ Etant donnés deux entiers $c\neq 0$ et $d$
de $\mc O_F$ premiers entre eux, nous définissons la somme de
Dedekind généralisée $s_j(d,c)$ comme une certaine fonction
réelle-analytique $s_j(d,c;\,.):\mc H^{n-1}\rightarrow \mb R.$

L'objectif de cet article est double : d'abord étudier ces
fonctions en généralisant les identités remarquables satisfaites
par les sommes de Dedekind classiques. Ensuite, relier des valeurs
spéciales de $s_j(d,c;\,.)$ en des points algébriques de $\mc
H^{n-1}$ à des valeurs spéciales de fonctions $L$ de Hecke en
$s=0.$

\smallskip

Commençons par présenter le cas où $F$ est le corps des
rationnels. Les sommes de Dedekind classiques peuvent être
introduites de la manière suivante à partir de la formule de
transformation modulaire du logarithme de la fonction $\eta$ de
Dedekind.

La fonction $\eta$ de Dedekind est définie sur le demi plan de
Poincaré $\mc H$ par $$\eta(z)=e^{i\pi z/12}\dis
\prod_{n=1}^\infty (1-e^{2i\pi nz}).$$ Sa puissance 24-ème est la
forme modulaire $\Delta$ de poids 12 sur le groupe $SL_2(\mb Z).$
On choisit comme logarithme de $\eta$ la branche holomorphe
\begin{equation*}\ln\eta(z):=\frac{i\pi}{12} z
-\sum_{m=1}^\infty\sum_{n=1}^\infty \frac{e^{2i\pi m nz}}
m.\end{equation*}
Soit $A=\left(\begin{array}{cc} a &b\\
c& d
\end{array}\right)$ une matrice de $ SL_2(\mb Z)$ telle que $c\neq0.$ Si on note $\ln$ la branche principale du logarithme, la fonction $\ln \eta$ vérifie donc la formule de transformation

\begin{equation*}\label{transfDedusu}
\ln\eta(Az)=\ln\eta(z)+ \frac 14 \ln\left(-(cz+d)^2\right) +\frac
{i\pi} {12} \Phi_R(A)
\end{equation*}
\noindent qui \textit{définit} la fonction $\Phi_R : SL_2(\mb
Z)\rightarrow \mb Z$ de Rademacher. La somme de Dedekind classique
$s(d,c)$ est le rationnel lié à $\Phi_R(A)$ par la relation
\begin{equation*}\label{PhiDEDusuel}
s(d,c):=\frac 1
{12}\left[-\sign(c)\Phi_R(A)+\frac{a+d}{|c|}\right]
.\end{equation*}

Dans [De], Dedekind déduit de la formule de transformation de $\ln
\eta$ la loi de réciprocité fondamentale
\begin{equation}\label{recipdedusu}
s(d,c)+s(c,d)=-\frac 14 +\dis\frac 1{12}\left(\frac d c+\frac cd
+\frac 1 {cd}\right)\ \ \ \ \ \textrm{si} \ \ (c,d)=1, \ \ \ c>0,\
d>0.
\end{equation}
Ensuite, il utilise le travail de Riemann pour démontrer l'égalité
\begin{equation*}\label{DEFDEDUSUSOMME}
s(d,c)=\sum_{k\, \textrm{mod}\, c} \left(\left(\frac k
c\right)\right)\left(\left( \frac {kd} c \right)\right),\ \
\textrm{ où } ((x))=\left\{
\begin{array}{ll}0 & \textrm{si}\ x \
\textrm{entier,}\\
x-[x]-\frac 12 &\textrm{sinon,}
\end{array}\right. \end{equation*}
qui permet d'étendre la définition de $s(d,c)$ au cas où $c$ et
$d$ ne sont pas premiers entres eux. Enfin Dedekind montre que
pour tout $p$ premier on a l'identité
\begin{equation}\label{DEDUSUPPRE}
s(dp,c)+\sum_{r\, \textrm{mod} \, p} s(d+cr,cp)=(p+1)s(d,c),
\end{equation}
c'est-à-dire que les sommes de Dedekind classiques sont des
fonctions propres pour certains opérateurs de Hecke.

Toute cette construction se généralise au corps $F.$ En effet,
nous utilisons la formule de transformation modulaire de la
fonction $\Lambda_j : \mc H^n\rightarrow \mb R$ introduite par
Hara ([Ha]) pour définir une fonction réelle-analytique $$\Phi_j :
SL_2(\mc O_F)\times \mc H^{n-1}\rightarrow \mb R$$ analogue de la
fonction $\Phi_R$ de Rademacher. Nous introduisons ensuite la
somme de Dedekind généralisée $s_j(d,c; \,.)$ qui ne diffère de
$\Phi_j \left(\left(\begin{array}{cc} a &b\\ c& d
\end{array}\right),\, .\right)$ que d'un terme élémentaire.

Nous démontrons alors que ces sommes vérifient une loi de
réciprocité et une identité qui généralisent (\ref{recipdedusu})
et (\ref{DEDUSUPPRE}). Cette loi de réciprocité, qui constitue le
théorème \ref{threcipDEDgene}, met en évidence le rôle particulier
joué par la fonction $s(0,1;\, .).$

\bigskip

On s'intéresse dans le reste de cet article à des invariants de
classes de $ \Gamma:=SL_2(\mc O_F)$ construits à l'aide de valeurs
spéciales de la fonction $\Phi_j.$

Cet invariant est défini dans le cas rationnel par Rademacher
([Ra1]) selon la formule
$\Psi(A):=\Phi_R(A)/6-\sign(c\,\tra(A))/2.$ Dans [At], M. Atiyah
identifie différents invariants de classes de $SL_2(\mb Z)$ à la
fonction $\Psi.$ En particulier, il explique comment associer à
une matrice hyperbolique de $SL_2(\mb Z)$ une fonction entière
$L_A(s)$ qui consiste essentiellement en une fonction $L$ de Hecke
partielle. M. Atiyah démontre ensuite l'égalité
$$L_A(0)=\Psi(A)$$ en adaptant des résultats de C. Meyer [Me].
Nous nous proposons de généraliser l'identité précédente en nous
inspirant des résultats de Hara [Ha].

On considère une matrice $A$ de $\Gamma$ comme un élément $(A_k)$
de $SL_2(\mb R)^n.$ On suppose que $A$ n'a qu'une composante
hyperbolique $A_j.$ Un telle matrice sera dite
\textit{quasi-elliptique}. En effet, les autres composantes sont
alors elliptiques et on note $\omega_c\in \mc H^{n-1}$ le points
fixe de $(A_k)_{k\neq j}.$

D'une part, on normalise la valeur spéciale $\Phi_j(A,\omega_c)$
pour obtenir un invariant $\Psi(A).$ On définit d'autre part la
fonction $L_A(s)$ associée à une matrice $A$ quasi-elliptique.
C'est pour l'essentiel une fonction $L$ de Hecke partielle d'un
ordre de l'extension quadratique $K$ de $F$ engendrée par les
valeurs propres de $A.$ Nous montrons ensuite que ces fonctions
$L_A(s)$ sont entières et qu'elles vérifient une équation
fonctionnelle en les reliant à la période de la dérivée partielle
d'une série d'Eisenstein réelle-analytique. En outre, ces
fonctions ont un zéro d'ordre $\geq n-1$ en $s=0.$

Nous concluons dans le théorème \ref{thL0PSI} à l'égalité
souhaitée
\begin{equation*}\label{LAPSIAresume}L_A^{(n-1)}(0)=(n-1)!\, \Psi(A)\end{equation*}
grâce à la formule limite de Kronecker généralisée.

En particulier, si $n=2,$ nous montrons que la nature arithmétique
de l'invariant $\Psi(A)$ est gouvernée par la conjecture de Stark
pour le corps $K$. On en déduit quelques valeurs de $\Psi(A).$

\smallskip

Nous expliquons enfin comment  ce travail permet d'interpr\'eter
la conjecture de Stark pour le corps $K$ comme compl\'ementaire
\`a la conjecture de Darmon [Da, Conjecture 8.17] pour les points
de Heegner de $K\cap\mc H.$ En effet, la construction de Darmon
est bas\'ee sur une forme modulaire de Hilbert $f$ cuspidale de
poids $(2,2)$. Cette construction conduit conjecturalement \`a des
points alg\'ebriques sur la courbe elliptique associ\'ee \`a $f.$

Notre article s'int\`egre dans le formalisme de Darmon en
remplaçant la forme cuspidale $f$ par la s\'erie d'Eisenstein de
poids $(2,2)$ pour $\Gamma.$ On est alors conduit \`a la valeur
sp\'eciale $L'_A(0)$ qui, si l'on en croit la conjecture de Stark,
est le logarithme d'une unit\'e alg\'ebrique.

\section{Sommes de Dedekind généralisées.}

Soit $F$ un corps de nombres totalement réel de nombre de classes
$1.$ On note $n=[F:\mb Q]$ son degré, $\mc O_F$ l'anneau des
entiers de $F,$ $\mc H$ le demi-plan de Poincaré. Fixons les $n$
plongements réels $\iota_1,\ldots,\iota_n$ de $F.$ Une matrice $A$
du groupe modulaire de Hilbert $\Gamma:=SL_2(\mc O_F)$ peut alors
être vue comme un élément $(A_k)$ de $SL_2(\mb R)^n.$ On en déduit
une action par homographies du groupe $\Gamma$ sur le produit $\mc
H^n.$

\smallskip

Nous étudions dans un premier temps une fonction définie sur $ \mc
H^n$ introduite par Hara dans [Ha, p.877]. Nous la considérons
comme un analogue pour $F$ du logarithme de la fonction $\eta$ de
Dedekind. La formule de transformation modulaire de cette fonction
nous permet ensuite de définir et d'étudier les sommes de Dedekind
généralisées associées à $F.$

\subsection{La fonction $\Lambda_j,$ analogue de $\ln
\eta.$}\label{DEFLAMBDAj}

\noindent On désigne par $d_F$ le discriminant de $F,$ $\mf
d=(\delta)$ la différente, $U_F$ (resp. $U_F^+$) le groupe des
unités (resp. unités totalement positives) de $F$ et $R_F$ son
régulateur. On notera $a_k:=\iota_k(a)$ l'image d'un élément $a$
de $F,$ et $x_k$ (resp. $y_k$) la partie réelle (resp. imaginaire)
d'un élément $z_k$ de $\mc H.$ En suivant [Ha], nous définissons
une fonction sur $\mc H^n$ qui sera l'analogue de $\ln \eta.$

\begin{defi}\label{defLAMBDAj}
Pour tout $j\in\{1,\ldots,n\}$, on définit la fonction $\Lambda_j$
en posant
\begin{equation*}\label{deflambdaj}
\begin{array}{cccl} \Lambda_j : & \mc H^n &\longrightarrow &\mb C\\
& z=(z_1,\ldots,z_n)& \longmapsto& \Lambda_j(z):= i\pi\kappa_F z_j
\left(\dis\prod_{k\neq j,\,k=1}^{n} y_k\right)-{\sqrt{d_F}\over
2R_F} \Omega_j(z),
\end{array} \end{equation*} où on a noté
$\kappa_F$ la constante $\frac{d_F\zeta_F(2)}{2^nR_F\pi^{n+1}},$
et $\Omega_j$ la fonction définie sur $\mc H^n$ par la série
absolument convergente
\begin{equation}\label{defOmegaj}\!\!\!\,\
\Omega_j(z):=\!\!\!\sum_{\nu \in
\mc{O}_F/{U_F^+}}\!\!\!\prim\!\!\!\! \frac{[U_F:
U_F^+]^{-1}}{|N_{F/\mb Q}(\nu)|} \!\!\!\sum_{ \scriptstyle{\mu \in
\mc{O}_F,}\atop \scriptstyle{\frac{\mu_j \nu_j}{
\mf{\delta}_j}>0}}\!\prim\!\!\!\!\!\!
e^{2i\pi\frac{\mu_j\nu_j}{\delta_j} z_j}\prod_{k\neq j,\,
k=1}^n\!\! e^{2i\pi\left(\frac{\mu_k\nu_k}{\delta_k}
x_k+i\left|\frac{\mu_k\nu_k}{\delta_k}\right|y_k\right)}.\!\end{equation}
Le symbole $\sum \ \prim\!\!\!\!\!\!$ signifie que la sommation
porte sur les éléments non nuls.
\end{defi}

\noindent Remarque : Dans cette somme, on a imposé la condition
$\frac{\mu_j \nu_j}{ \mf{\delta}_j}>0.$ On a donc construit
$n=[F:\mb Q]$ fonctions $\Lambda_j$ en privilégiant successivement
chaque plongement $\iota_j$ de $F.$ Une conséquence immédiate de
ce choix est que \textit{$\Lambda_j(z_1,\ldots,z_n)$ est
holomorphe par rapport à $z_j$, mais pas par rapport à $z_k$ pour
$k\neq j.$}

\smallskip

En fait, $\Lambda_j(z)$ apparaît naturellement comme le morceau
holomorphe par rapport à $z_j$ de la partie réelle de
$\Lambda_j(z).$ Cette fonction $\Reel \Lambda_j,$ qui ne dépend
pas de l'entier $j,$ a déjà été étudiée par Asai dans \cite{As}.
Il démontre qu'elle est l'analogue pour $F$ de $\ln|\eta|=\Reel
(\ln \eta).$ Ceci nous conduit à considérer que $\Lambda_j$ est
l'analogue pour $F$ de $\ln \eta.$

Rappelons plus précisément les résultats obtenus par Asai. En
suivant ses notations, on pose
\begin{equation}\label{lienhLambda}h(z):=-4\Reel\left(\Lambda_j(z)\right).\end{equation}
\noindent On doit aussi introduire la série d'Eisenstein
non-holomorphe associée à $F.$

\begin{defi}\label{defEisnonholo}
Soit $z=(z_1,\ldots,z_n) \in \mc H^n.$ Pour $\Reel(s)>1$ on
définit la \textit{série d'Eisenstein non-holomorphe} $E_F(z,s)$
par la série absolument convergente
$$E_F(z,s):=\sum_{(\mu,\nu)\in \mc O_F^2/U_F}\!\!\!\!\!\!\prim
\prod_{k=1}^n \frac{ y_k^s }{|\mu_kz_k+\nu_k|^{2s}}.$$
\end{defi}

Les propriétés fondamentales de la série $E_F(z,s)$ et la
``formule limite'' qui la relie à $h(z)$ ont déjà été établies
dans [As, p. 204 et Th. 3]. Ces résultats classiques sont
rassemblés dans le théorème qui suit.

\smallskip
\begin{theo} [Asai]~
\label{propEasai}
\begin{enumerate}
\item Définie par une série absolument convergente pour
$\Reel(s)>1,$ la fonction $s\mapsto E_F(z,s)$ se prolonge en une
fonction holomorphe sur tout le plan complexe sauf en $s=1$ où
elle a un pôle simple. \item Pour toute matrice $A$ de $\Gamma,$
on a $E_F(Az,s)=E_F(z,s).$ \item Pour tout $ z\in \mc H^n,$
$E_F(z,s)$ vérifie l'équation fonctionnelle
\begin{equation*}G_F(2s)E_F(z,s)=G_F(2-2s)E_F(z,1-s),\end{equation*}
où $G_F(s):= d_F^{\frac s2}\pi^{-\frac {ns} 2} \Gamma(\frac s2)^n$
est le facteur Gamma de la fonction zeta de Dedekind du corps $F.$
\item ``Formule limite de Kronecker généralisée.''

Au voisinage de $s=1,$ on a le développement en série de Laurent
\begin{equation*}
E_F(z,s)=\frac{(2\pi)^nR_F}{4d_F}\left[\frac 1 {s-1}+ \gamma_F-
\ln\left(\prod_{k=1}^n
y_k\right)+h(z)\right]+O(s-1),\end{equation*} où $\gamma_F$ est la
constante $\frac{4\sqrt{d_F}}{2^\frac n2
R_F}((s-1)\zeta_F(s))'(1)-\ln 2^n.$

\end{enumerate}
\end{theo}\noindent
Récapitulons aussi les propriétés de la fonction $h$ données dans
[As, Th. 4-5].

\smallskip

\begin{theo}[Asai] ~\label{propReLambda}
\begin{enumerate}
\item La fonction $h(z)$ est une fonction pluri-harmonique à
valeurs réelles.

\item Pour toute matrice $A=\left(\begin{array}{cc} a &b\\ c& d
\end{array}\right)\in \Gamma,$ on a la relation modulaire
$$h(z)=h(Az)+\sum_{k=1}^n
\ln|c_kz_k+d_k|^2.$$

\item Les fonctions $h(z)$ et $\zeta_F(s)\zeta_F(s+1)$ sont
associées via la transformée de Mellin.

\item On a  l'\'egalit\'e $h(z)=-4\ln |\eta(z)|$ dans le cas o\`u $F$ est le corps des rationnels.
\end{enumerate}
\end{theo}

Ces quatre propriétés ainsi que la formule limite de Kronecker
généralisée permettent à T. Asai de considérer que
$-h/4=\Reel(\Lambda_j)$ est l'analogue pour $F$ de la fonction
$\ln|\mathbf \eta|=\Reel\left(\ln \eta\right).$

\medskip

Sous l'action du groupe modulaire de Hilbert, la partie réelle de
$\Lambda_j(z)$ vérifie une formule explicite donnée par le
théorème \ref{propReLambda}.2. Nous allons en déduire que la
partie imaginaire de $\Lambda_j$ vérifie une formule de
transformation où apparaît naturellement une fonction $\Phi_j$ à
valeurs réelles.

\smallskip

Pour cela, on considère la fonction $\boldsymbol{\eta}_{j}(z):=
e^{ \Lambda_j(z)}$ analogue de la fonction $\eta$ de Dedekind. Par
construction, elle est holomorphe par rapport à $z_j$ et elle ne
s'annule pas sur $\mc H^n.$

\noindent Etudions son comportement
sous l'action d'une matrice $ A= \left(\begin{array}{cc} a &b\\
c& d
\end{array}\right)$ de $\Gamma.$

\smallskip

Tout d'abord, on remarque que si $c=0$ la formule de
transformation de $\boldsymbol{\eta}_j(z)$ sous l'action de $A$
est complètement explicite. En effet, l'entier $a=d^{-1}$ est
alors une unité de $F.$ La transformation modulaire associée est
donc du type $z\mapsto a^2z+a b.$ Ces transformations laissent la
fonction $\Omega_j(z)$ définie par (\ref{defOmegaj}) invariante.
On déduit de la définition \ref{deflambdaj} que
$$\boldsymbol{\eta}_j^4\left(\frac{az+b} d \right)=\boldsymbol{\eta}_j^4(z)e^{4i\pi\kappa_Fb_jd_j\prod_{k\neq j} y_k}.$$

\smallskip

Revenons au cas général. Pour une matrice $A$ quelconque du groupe
$\Gamma,$ on sait d'après le théorème \ref{propReLambda}.2 que la
fonction $| \boldsymbol{\eta}_{j}|=e^{\Reel(\Lambda_j)}$ vérifie
\begin{equation}\label{modetanew}
|\boldsymbol{\eta}_{j}(Az)|^4=|
\boldsymbol{\eta}_{j}(z)|^4|c_jz_j+d_j|^2\prod_{k\neq j,k=1
}^n|c_kz_k+d_k|^2.\end{equation} Le lemme suivant résulte
immédiatement de cette formule.
\begin{lem} \label{lemphidef} Soient $A=
\left(\begin{array}{cc} a &b\\ c& d
\end{array}\right)$ une matrice de $ \Gamma $ et $z=(z_1,\ldots,z_n)\in\mc H^n.$ Il existe un réel $\phi$
indépendant de $z_j \in \mc H$ tel que

$$\boldsymbol{\eta}_{j}^4 (Az)= (c_jz_j+d_j)^2\left[
\prod_{k\neq j,k=1
}^n|c_kz_k+d_k|^2\right]\boldsymbol\eta_{j}^4(z)e^{4i\pi
\phi+i\pi}.$$
\end{lem}

\begin{proof}

Le quotient $$\dis\frac{\dis \boldsymbol{\eta}_j^4( Az)}{\dis
(c_jz_j+d_j)^2\boldsymbol{\eta}_j^4(z) \prod_{k\neq j,k=1}^n
|c_kz_k+d_k|^2}$$ est une fonction de la variable $z_j$ holomorphe
sur le connexe $\mc H.$ Son module est égal à 1 d'après
(\ref{modetanew}). Ce quotient est donc égal à une constante (par
rapport à $z_j$) de module 1.
\end{proof}

\smallskip

Le réel $\phi$ du lemme précédent ne dépend donc que de la matrice
$A$ et de $\hat{z}_j:=(z_1,\ldots,z_{j-1},z_{j+1},\ldots,z_n)\in
\mc H^{n-1}.$ Ce réel étant pour l'instant défini à un demi-entier
près, nous allons maintenant le définir de manière univoque en
choisissant bien sûr $\Lambda_j$ comme logarithme de
$\boldsymbol{\eta}_j.$

\begin{defi} \label{defiPHI} Soit $z=(z_1,\ldots,z_n)\in \mc H^n,$ et $\hat{z}_j\in \mc H^{n-1}$ le $(n-1)$-uplet
$\hat{z}_j:=(z_1,\ldots,z_{j-1},z_{j+1},\ldots,z_n)$ associé. Pour
tout $j\in\{1,\ldots,n\},$ on définit la fonction
$$\begin{array}{ccccll} \Phi_j :& \Gamma
& \times& \mc H^{n-1}&\longrightarrow & \mb R\\
& A=\left(\dis\begin{array}{cc} a &b\\ c& d
\end{array}\right)&,&\hat{z}_j
&\longmapsto& \Phi_j(A,\hat{z}_j),\end{array}$$
$$\textrm{où } \ \ \ \Phi_j(A,\hat{z}_j):= \left\{
\begin{array}{ll}\kappa_F b_jd_j\dis\prod_{k\neq
j} y_k \ \ \textrm{ si } c=0,& \\
\dis\frac{1}{\pi}\Ima \! \left[\Lambda_j(Az)-\Lambda_j(z)-\frac
14\ln\left[-(c_jz_j+d_j)^2\right]\right] & \,\textrm{si } c\neq
0.\end{array}\right.
$$

\end{defi}\noindent
Le lemme \ref{lemphidef} assure que cette définition a un sens.

\medskip
\noindent La fonction $\Phi_j(A,\hat{z}_j)$ est une généralisation
de la \textit{fonction $\Phi_R : SL_2(\mb Z)\rightarrow \mb Z$ de
Rademacher} ([Ra2 p.150]) qui apparaît naturellement dans la
formule de transformation du logarithme de la fonction $\eta$ de
Dedekind. Notons que notre normalisation donne dans le cas
rationnel $\Phi_1=-\Phi_R/12.$

\medskip

En conclusion, la fonction $\Lambda_j$ se transforme sous l'action
d'une matrice quelconque $ A=\left(\dis\begin{array}{cc} a &b\\ c&
d
\end{array}\right)$ de $\Gamma$ selon la règle
\begin{equation}\label{LAMBDATRANSFPHI}
\Lambda_j(Az)=\Lambda_j(z)+\frac{\delta_c} 4
\left[\ln\left[-(c_jz_j+d_j)^2\right]+\sum_{k\neq j}\ln
|c_kz_k+d_k|^2\right]+i\pi\Phi_j(A,\hat{z}_j),\end{equation} où
$\delta_c$ vaut $0$ si $c=0$ et $1$ sinon.
\medskip




\bigskip

Etant donnés une matrice $A=\left(\begin{array}{cc} * &*\\
c& d
\end{array}\right)$ de $SL_2(\mb R)$ et $z\in \mc H,$ on note $L(A,z)$ la fonction $L(A,z) := \delta_c\Ima \ln \left[-(cz+d)^2\right].$
On pose alors $$-\pi\varDelta(A,B):=L(AB,z)-L(A,Bz)-L(B,z).$$ La
fonction $\varDelta$ ainsi définie ne dépend pas du choix de $z\in
\mc H.$ En outre, $\varDelta$ est un $2$-cocycle sur $SL_2(\mb R)$
à valeurs dans $\mb Z.$ Ces deux assertions se déduisent
immédiatement du lemme qui suit.
\smallskip
\begin{lem}\label{lemmedelta} Soit $A=\left(\begin{array}{cc} * &*\\
c& d
\end{array}\right)$ et $B=\left(\begin{array}{cc} * &*\\ c'& d'
\end{array}\right)$ deux matrices de $SL_2(\mb R).$ On note $AB=\left(\begin{array}{cc} * &*\\ c''&
d''
\end{array}\right)$ leur produit.
On a alors l'égalité

\begin{equation}\label{butth1}\varDelta(A,B)=-\sign(cc'c''),\end{equation}

où $\sign(x)$ vaut bien sûr $0$ si $x=0,$ $1$ si $x>0$ et $-1$ si
$x<0.$
\end{lem}

\begin{proof}
\noindent Un calcul direct permet d'obtenir l'identité
\begin{equation}\label{liencc'c''}
(cBz+d)(c'z+d')=(c''z+d'').
\end{equation}

\noindent Nous allons exprimer cette égalité en termes de
logarithmes. La démonstration du lemme repose sur les propriétés
élémentaires de la branche principale du logarithme :

i) $\ln(zz')=\ln( z)+\ln (z')$ à condition que $|\!\textrm{ Arg
}(z)+\!\textrm{ Arg }(z')|<\pi.$

ii) $\ln (z^{-1})=-\ln (z).$

\noindent Nous devons distinguer trois cas.
\begin{enumerate}
\item
Commençons par supposer que deux des trois réels $c,c',c''$ sont
nuls. Dans ce cas, ils sont tous les trois nuls car les matrices
triangulaires supérieures forment un groupe. Ainsi
$\delta_{c}=\delta_{ c'}=\delta_{ c''}=0$ et les deux membres de
(\ref{butth1}) sont nuls donc égaux.

\item On traite à présent le cas où $c=0$ et $c'c''\neq 0.$
L'égalité (\ref{liencc'c''}) se réduit alors à
$d(c'z+d')=(c''z+d'').$ On déduit donc de i) que
$$\ln(d^2)+\ln\left[-(c'z+d')^2\right]=\ln \left[-(c''z+d'')^2\right].$$
En prenant la partie imaginaire, on en conclut que l'égalité
(\ref{butth1}) est vérifiée. Le cas où $c'=0$ et le cas où $c''=0$
se traitent de façon analogue en utilisant les propriétés i) et
ii).

\item Il reste à prouver (\ref{butth1}) dans le cas où $cc'c''\neq
0.$

On note d'abord que les deux membres de l'égalité souhaitée ne
dépendent que de la classe de $A$ et $B$ dans $PSL_2(\mb R).$ Par
suite, quitte à changer $A$ ou $B$ en leur opposé, on peut
supposer que $c>0$ et $c'>0.$

Sous cette hypothèse, les complexes $-i(cBz+d)$ et $ -i(c'z+d') $
sont deux éléments du demi-plan $\Reel(z)>0,$ donc leur argument
est dans l'intervalle $]-\frac\pi 2,\frac \pi 2[.$ Ainsi, on
déduit de i) et (\ref{liencc'c''}) que
$$L(A,Bz)+L(B,z)=2\Ima \ln\left[-(c''z+d'')\right].$$
La définition de $\varDelta$ permet de conclure que
\begin{equation}\label{lastRln}\pi\varDelta(A,B)=-L(AB,z)+2\Ima
\ln\left[i^2(c''z+d'')\right].\end{equation} Distinguons pour
finir les deux sous-cas $c''<0$ et $c''>0.$
\begin{itemize}
\item Sous-cas 3.1 : si $c''<0,$ alors le complexe $ i(c''z+d'')
$ est dans le demi-plan $\Reel(z)>0.$ On a donc d'après la
propriété i) : $$L(AB,z)=2\Ima\ln (i(c''z+d''))=2\Ima\left(\ln
\left[i^2(c''z+d'')\right]-\ln (i)\right).$$ Par suite, l'égalité
(\ref{lastRln}) devient
$$\pi\varDelta(A,B)=2\Ima \ln (i)=\pi .$$ \item
Sous-cas 3.2 : si $c''>0,$ alors $\Reel(-i(c''z+d''))>0. $ On en
conclut que $\mc \pi\varDelta(A,B)=2\Ima \ln (-i)=-\pi$ de façon
similaire au sous-cas précédent.
\end{itemize}

\end{enumerate}

On réunit ces deux sous-cas en écrivant
$\varDelta(A,B)=-\sign(cc'c'').$
\end{proof}

\medskip

Le $2$-cocycle $\varDelta$ a une interprétation géométrique. (Je
remercie E. Ghys de m'avoir signalé ce fait). Pour voir ceci,
identifions $\mc H$ avec le disque de Poincaré et fixons
arbitrairement un point $x$ sur le cercle à l'infini. Le réel $\pi
\varDelta(A,B)$ apparaît alors comme l'aire algébrique du triangle
idéal de sommets $x,$ $Ax,$ $ABx.$ Par suite, le cocycle
$\varDelta$ est appelé le \textit{$2$-cocyle d'aire} de $SL_2(\mb
R)$ (voir [K-M, p. 238]).

\smallskip
En composant chacun des $n$ plongements de $\Gamma$ dans $SL_2(\mb
R)$ avec le cocycle d'aire $\varDelta,$ on obtient donc $n$
cocycles distincts sur le groupe $\Gamma.$

Notons en particulier $\varDelta_{SL_2(\mb Z)}$ la restriction de
$\varDelta$ à $SL_2(\mb Z).$ Sa classe de cohomologie est un
élément d'ordre 3 de $H^2(SL_2(\mb Z),\mb Z)=\mb Z/12\mb Z.$ On
sait de plus que le premier et le deuxième groupe de cohomologie
rationnelle $H^1(SL_2(\mb Z),\mb Q)$ et $H^2(SL_2(\mb Z),\mb Q)$
sont nuls. Il existe donc un unique $1$-cocycle \textit{à valeurs
rationnelles} sur $SL_2(\mb Z)$ dont le cobord est
$\varDelta_{SL_2(\mb Z)}.$ Ce $1$-cocycle n'est autre que
$-\Phi_R/3,$ comme le montre le résultat de Rademacher ([Ra2
p.152]) : pour toutes matrices $A$ et $B$ de $SL_2(\mb Z),$ on a
$$\Phi_R(AB)-\Phi_R(A)-\Phi_R(B)=3\sign(cc'c'').$$

Nous proposons une généralisation de cette identité fondamentale
dans le théorème suivant.
\bigskip

\begin{theo}\label{phicomp}
Soit $n=[F:\mb Q]$ et $j\in\{1,\ldots,n\}.$ Pour tout $\hat{z}_j
\in \mc H^{n-1}$ et toutes matrices $A$ et $B$ de $ \Gamma$, on a
la relation
\begin{equation*} \label{PhiAB} \Phi_j(AB,\hat{z}_j)-\Phi_j(
A,\widehat{Bz}_j)-\Phi_j(B,\hat{z}_j)=-\frac 14
\sign(c_jc_j'c_j''),\end{equation*}

où on a noté $\widehat{Bz}_j$ le $(n-1)$-uplet
$(B_1z_1,\ldots,B_{j-1}z_{j-1},B_{j+1}z_{j+1},\ldots,B_nz_n).$

\end{theo}

\smallskip

\noindent \textit{Démonstration du théorème \ref{phicomp}.}
L'égalité (\ref{LAMBDATRANSFPHI}) permet d'écrire successivement
\begin{equation*}\label{lambdaAB1} \Ima \Lambda_j((AB)
z)=\Ima (\Lambda_j(z))+\frac {1}4L(A_jB_j,z_j)+\pi
\Phi_j(AB,\hat{z}_j),\end{equation*}
\begin{equation*} \label{LambdaA1B} \Ima \Lambda_j(A(
Bz))=\Ima (\Lambda_j(B z))+\frac {1} 4L(A_j,B_jz_j)+\pi \Phi_j(
A,\widehat{Bz}_j),\end{equation*}
\begin{equation*}\label{LAMBDABSEUL} \Ima \Lambda_j(B z)=\Ima(\Lambda_j(z))+\frac {1}4L(B_j,z_j)
+\pi\Phi_j(B,\hat{z}_j).
\end{equation*}
En combinant ces trois équations, il vient
\begin{equation*}\label{phiABfction log}4\pi\left( \Phi_j(AB,\hat{z}_j)-\Phi_j(A,\widehat{Bz}_j)-\Phi_j(
B,\hat{z}_j)\right)=\pi\varDelta(A_j,B_j).\end{equation*} Le
résultat souhaité se déduit donc du lemme \ref{lemmedelta}. \hfill
$\Box$

\medskip

Remarque : fixons un point $\omega_c\in \mc H^{n-1},$ et
consid\'erons le sous-groupe $\Gamma_{\omega_c}$ de $\Gamma$
constitu\'e des matrices $A$ telles que $\omega_c$ est un point
fixe de $(A_k)_{k\neq j}.$ En choisissant $\hat{z}_j=\omega_c$
dans le th\'eor\`eme pr\'ec\'edent, on voit que la restriction du
$\varDelta$ \`a $\iota_j(\Gamma_{\omega_c})$ est le cobord du
1-cocycle $\Phi_j(.,\omega_c).$ Ce sous-groupe $\Gamma_{\omega_c}$
est tr\`es petit : on verra dans le paragraphe \ref{fctionL} qu'il
est de rang $leq 1$. Le cocycle $\Phi_j(.,\omega_c)$ de
$\Gamma_{\omega_c}$ sera d'une grande importance dans la partie
\ref{partieL} : nous le relierons \`a la valeur sp\'eciale de
certaines fonctions $L$ de Hecke en $s=0$ chaque fois que
$\Gamma_{\omega_c}$ est de rang 1.

\subsection{Sommes de Dedekind
généralisées.}\label{DEDgene} La formule de transformation de la
fonction $\Lambda_j$ va nous permettre de définir formellement les
sommes de Dedekind associées au corps $F.$ Notre construction
procède de manière analogue au cas du logarithme de la fonction
$\eta$ et des sommes de Dedekind classiques. Etant donnés deux
entiers $c\neq 0$ et $d$ de $\mc O_F$ premiers entre eux, nous
définissons la somme de Dedekind généralisée associée comme une
fonction $s(d,c;\,.) : \mc H^{n-1}\rightarrow \mb R.$

Nous montrons ensuite que les sommes généralisées vérifient une
loi de réciprocité. Cette loi, qui fait l'objet du théorème
\ref{threcipDEDgene}, met en évidence la somme particulière
$s(0,1;z_2,\ldots,z_n).$ Plus précisément, la loi de réciprocité
nous permet d'exprimer toute somme de Dedekind généralisée comme
une somme finie de $s(0,1;.)$ et de termes élémentaires.

Dans la proposition \ref{s01prop}, nous utilisons le travail de
Hecke [He] pour écrire la fonction fondamentale $s(0,1;z_2)$ sous
la forme d'une série remarquable. Le terme général de cette série
fait intervenir des valeurs spéciales de fonctions $L$ de Hecke
sur la droite $\Reel (s)=1.$

Enfin, nous montrons que les sommes de Dedekind généralisées sont
vecteurs propres de certains opérateurs de Hecke.

\subsubsection{Définition et premières propriétés.}\label{pardefded}
Nous avons introduit la fonction $\Phi_j(A,\hat{z}_j)$ qui
apparaît naturellement dans la formule de transformation modulaire
de $\Lambda_j(z).$ Comme dans le cas rationnel, nous la
normalisons à l'aide d'un facteur élémentaire pour définir les
sommes de Dedekind généralisées.

\begin{defi} \label{definitionsommesgene} Soit $c\neq 0 $ et $d$ des entiers de $\mc O_F$ premiers entre eux.
Pour tout $j\in \{1,\ldots,n\}$, la \textit{somme de Dedekind
généralisée $s_j$ associée à $(c,d)$} est la fonction $s_j(d,c;.)
: \mc H^{n-1} \rightarrow \mb R$ définie par la formule

\begin{equation} \label{defsommesdedgene}
s_j(d,c;\hat{z}_j) =-\sign(c_j)\Phi_j\left(\left(\begin{array}{cc} a &b\\
c& d
\end{array}\right),\hat{z}_j\right)+
\frac{
\kappa_F}{|c_j|}\left[a_jf_j(d,c;\hat{z}_j)+d_jf_j(0,1;\hat{z_j})\right],
\end{equation}

\noindent où $\kappa_F=\frac{d_F\zeta_F(2)}{2^nR_F\pi^{n+1}},$ $a$
et $b$ étant des entiers de $\mc O_F$ tels que $ad-bc=1.$ On a
noté $f_j(d,c;\hat{z}_j)$ la fonction élémentaire
$$f_j(d,c;\hat{z}_j)=\prod_{k\neq j} \frac{y_k}{|c_kz_k+d_k|^2}.$$
\end{defi}

Pour que $s_j(d,c;\hat{z}_j)$ soit bien défini, nous devons
maintenant vérifier que le membre de droite de l'égalité
(\ref{defsommesdedgene}) ne dépend pas du choix des entiers $a$ et
$b$ de $\mc O_F$ qui vérifient $ad-bc=1.$

\noindent Notation : quitte à renuméroter les plongements, on
suppose désormais que $j=1;$ on note $s$ la somme $s_1$ et
$\Lambda,$ $\Omega,$ $\Phi,$ $f$ les fonctions $\Lambda_1,$ $
\Omega_1,$ $\Phi_1$ et $f_1.$

\medskip

\noindent La méthode de Riemann-Dedekind ([De]) suggère de poser
$z_1=-\frac{d_1}{c_1}+i\frac t {|c_1|},$ où $t$ est un réel
positif qui va tendre vers 0. En notant $A$ la matrice $\left(\begin{array}{cc} a&b\\
c& d
\end{array}\right)$ de $\Gamma,$ la formule de transformation
(\ref{LAMBDATRANSFPHI}) donne l'égalité
\begin{equation}\label{lambdatend0}\!
\Ima\!\left[\Lambda\left(\frac{a_1}{c_1}+\frac{i}{t|c_1|},\widehat{Az}_1\right)\right]=
\Ima\!\left[\Lambda\left(-\frac{d_1}{c_1}+i\frac{t}{|c_1|},\hat{z}_1\right)\right]+\pi\Phi\left(A,\hat{z}_1\right).
\end{equation}
\noindent Rappelons que $\Lambda(z)=i\pi\kappa_F z_1
f(0,1;\hat{z_1})-\frac{\sqrt{d_F}}{2R_F} \Omega(z),$ où $\Omega$
est donnée par (\ref{defOmegaj}).

Lorsque $t$ tend vers $0,$ le complexe
$$\Omega\left(\frac{a_1}{c_1}+\frac{i}{t|c_1|},\frac{a_2z_2+b_2}{c_2z_2+d_2},\ldots,\frac{a_nz_n+b_n}{c_nz_n+d_n}\right)$$
tend vers $0.$ Ainsi le membre de gauche de l'égalité
(\ref{lambdatend0}) a pour limite
$$\frac{a_1} {c_1}\pi \kappa_F f(d,c;\hat{z}_1).$$ Quant au membre de droite, il a pour
limite
$$-\frac{d_1} {c_1}\pi \kappa_F
f(0,1;\hat{z}_1)-\frac{\sqrt{d_F}}{2R_F}\lim_{t\rightarrow 0\atop
t>0}
\Ima\left[\Omega\left(-\frac{d_1}{c_1}+it,\hat{z}_1\right)\right]+
\pi\Phi\left(A,\hat{z}_1\right).$$ On rassemble ces égalités pour
écrire finalement :
$$\frac{\kappa_F} {c_1} \left[a_1 f(d,c;\hat{z}_1) +d_1f(0,1;\hat{z}_1)\right] =
\Phi\left(A,\hat{z}_1\right)-\frac{\sqrt{d_F}}{2\pi
R_F}\lim_{t\rightarrow 0\atop t>0}
\Ima\,\Omega\left(-\frac{d_1}{c_1}+it,\hat{z}_1\right).$$

On conclut de cette égalité que le membre de droite de
(\ref{defsommesdedgene}) ne dépend pas du choix de $a$ et $b.$ Il
s'ensuit que $s_j(d,c; \hat{z}_j)$ est bien défini. De plus on a
obtenu la proposition suivante :
\smallskip
\begin{prop}\label{propdefsommes} Soit $c\neq 0$ et $d$ deux entiers de $\mc O_F$ premiers entre
eux, et $\hat{z}_1=(z_2,\ldots,z_n)\in \mc H^{n-1}.$ On a alors
\begin{equation*}
s(d,c;\hat{z}_1)=-\frac{\sign(c_1)\sqrt{d_F}}{2\pi
R_F}\lim_{t\rightarrow 0\atop t>0}
\Ima\,\Omega\left(-\frac{d_1}{c_1}+it,\hat{z}_1\right).
\end{equation*}
Cette égalité permet de définir $s(d,c,\hat{z_1})$ même si $c$ et
$d$ ne sont pas premiers entre eux. On a alors pour tout entier
$\lambda\neq 0$ de $\mc O_F$ :
$$s( \lambda d,\lambda c;\hat{z}_1)=s(d,c;\hat{z}_1).$$
\end{prop}

On peut donner tout de suite quelques identités élémentaires pour
les sommes généralisées. Elles sont tout à fait analogues au cas
des sommes de Dedekind classiques, si ce n'est la dépendance en
$z_2,\ldots,z_n$.

\smallskip
\begin{prop}\label{Dsigne}~ Soient $c\neq 0$ et $d$ des entiers de $\mc O_F.$ On a pour tout
$\hat{z}_1=(z_2,\ldots,z_n)\in \mc H^{n-1}$ les égalités :
\begin{enumerate}
\item
\begin{equation}\label{sde-c}
s(d,-c;\hat{z}_1)=s(d,c;-\overline{\hat{z}_1}).\end{equation}
\item
\begin{equation}\label{CONJ-D}s(-d,c;\hat{z}_1)=-s(d,c;-\overline{\hat{z}_1}).\end{equation}

\item Soit $\epsilon$ une unité de $F$ telle que
$\epsilon_1>0.$ Etant donné $z_k=x_k+iy_k\in \mc H,$ on désigne
par $|\epsilon_k|.z_k\in \mc H$ le complexe
$\epsilon_kx_k+i|\epsilon_k|y_k.$ Alors on a
\begin{equation}\label{sdcepsi} s(d,\epsilon c;\hat{z}_1)
=s(d,c;|\epsilon_2|.z_2,\ldots,|\epsilon_n|.z_n).\end{equation}

\item Pour tout entier $q$ de $\mc O_F,$ on a
\begin{equation}s(d+q c,c;z_{2}+q_2,\ldots,z_n+q_n)=s(d,c;z_{2},\ldots,z_n).\label{stranslation}\end{equation}

\end{enumerate}
\end{prop}

\begin{proof}
D'une part, il résulte immédiatement de la proposition
\ref{propdefsommes} que $s(-d,-c;\hat{z}_1)=s(d,c;\hat{z}_1).$ On
a d'autre part l'identité facile
$\Omega(-\bar{z})=\overline{\Omega(z)}.$ Les deux premières
assertions s'ensuivent à l'aide de la proposition
\ref{propdefsommes}.

\smallskip

Pour démontrer l'assertion $3,$ commençons par remarquer qu'en
changeant $\mu$ en $\epsilon\mu$ dans (\ref{defOmegaj}), on
obtient
$\Omega(|\epsilon_1|.z_1,\ldots,|\epsilon_n|.z_n)=\Omega(z).$ La
proposition \ref{propdefsommes} permet alors de conclure.
\noindent On procède de même à partir de l'égalité facile
$\Omega(z+q)=\Omega(z)$ pour établir l'assertion 4.
\end{proof}

\subsubsection{Loi de réciprocité des sommes de Dedekind généralisées.}

Nous montrons que les sommes de Dedekind généralisées introduites
précédemment vérifient une loi de réciprocité. Ensuite, nous
expliquerons comment cette loi réduit l'étude d'une somme
quelconque $s(d,c;\,\hat{z}_1)$ à celle de la fonction $s(0,1;\,
.):\mc H^{n-1}\rightarrow \mb R.$

\medskip

\begin{theo}[Loi de réciprocité]~\label{threcipDEDgene}
Soit $(c,d)$ un couple d'entiers de $\mc O_F$ premiers entre eux
tels que $c_1>0$ et $d_1>0.$ Notons
$\kappa_F:=\frac{d_F\zeta_F(2)}{2^nR_F\pi^{n+1}}.$ On a pour tout
$\hat{z}_1=(z_2,\ldots,z_n)\in \mc H^{n-1}$ l'identité

\begin{align*}s(d,c;\hat{z}_1)+s(c,d;\overline{\hat{z}_1^{-1}})
= s(0,1;\hat{z}_1)-\frac 1 4 &+\kappa_F \left[\frac{\dis d_1}{\dis
c_1}+\frac{\dis c_1}{\dis
d_1}\dis\prod_{k=2}^n\left|z_k\right|^{-2} \right.\\& \left.
\qquad \ +\frac{\dis 1}{\dis
c_1d_1}\prod_{k=2}^n\left|c_kz_k+d_k\right|^{-2}\right]\prod_{k=2}^ny_k.
\end{align*}

\end{theo}

\medskip\noindent
Remarque : Ce théorème généralise la loi de réciprocité
(\ref{recipdedusu}) des sommes de Dedekind classiques. On rappelle
que si $F=\mb Q$ on a $s(0,1)=0$ et $\kappa_\mb Q=1/12.$

\smallskip
\noindent \textit{Démonstration du théorème \ref{threcipDEDgene}.}
Soient $c$ et $d$ premiers entre eux tels que $c_1>0$ et $d_1>0.$
Soit $(a,b)\in \mc O_F\times \mc O_F$ vérifiant $ad-bc=1.$ La loi
de réciprocité se déduit du théorème \ref{phicomp} avec un choix
judicieux de matrices $A$ et $B.$

Choisissons en effet $A=\left(\begin{array}{cc} a &b\\
c& d
\end{array}\right),$ $B=\left(\begin{array}{cc} 0 &-1\\ 1& 0
\end{array}\right)$ dans ce théorème. On en déduit une identité qui s'écrit en termes de sommes de Dedekind
généralisées sous la forme\begin{align*}\label{phicompprloirecip}
s(-c,d;\hat{z}_1)-s(d,c;-\hat{z}_1^{-1})= s(0,1;\hat{z}_1)+\frac 1
4 &-\kappa_F
\left[\frac{c_1}{d_1}+\frac{d_1}{c_1}\prod_{k=2}^n|z_k|^{-2}\right.\\&\left.\qquad
+\frac{\prod_{k=2}^n
|c_k-d_kz_k|^{-2}}{c_1d_1}\right]\prod_{k=2}^ny_k.\end{align*}

\noindent On veut changer $\hat{z}_1$ en $-\hat{z}_1^{-1}$ dans
cette égalité. En utilisant à nouveau le théorème \ref{phicomp}
avec le couple de matrices $(B,B^{-1}),$ on obtient déjà
\begin{equation*} s(0,1;-\hat{z}_1^{-1})=-s(0,1;\hat{z}_1).
\end{equation*} Il résulte des deux identités précédentes que
\begin{align*}
s(d,c;\hat{z}_1)-s(-c,d;-\hat{z}_1^{-1})=s(0,1;\hat{z}_1)-\frac 1
4&+\kappa_F \left[\frac{c_1}{\dis
d_1}\prod_{k=2}^n|z_k|^{-2}+\frac{d_1}{c_1}\right.\\ &\left. \ \
\, +\frac{1}{\dis c_1d_1}\prod_{k=2}^n
|c_kz_k+d_k|^{-2}\right]\prod_{k=2}^ny_k.\end{align*} \noindent On
conclut à l'identité souhaitée grâce à la relation
(\ref{CONJ-D}).\hfill $\Box$

\subsubsection{La somme fondamentale
$s(0,1;z_2,\ldots,z_n).$}\label{sommededfondpart} La loi de
réciprocité montre que la fonction $s(0,1;\hat{z}_1)$ définie sur
$\mc H^{n-1}$ joue un rôle privilégié. On se propose dans ce
paragraphe de l'étudier plus en détail. La proposition suivante
justifie le nom de \textit{somme de Dedekind généralisée
fondamentale} pour cette fonction.

\medskip
\begin{prop}\label{remeucl} Soit $F$ totalement réel de nombre de classes
$1.$ La loi de réciprocité et la proposition \ref{Dsigne}
permettent d'exprimer toute somme de Dedekind généralisée
$s(d,c;\hat{z}_1)$ comme une somme finie de sommes $s(0,1;.)$ et
de termes élémentaires.
\end{prop}

\begin{proof}

Commençons par supposer que $\mc O_F$ est un anneau euclidien pour
la norme. La preuve consiste en un algorithme qui est calqué sur
\textit{l'algorithme d'Euclide.}

\smallskip
La boucle principale de l'algorithme est la suivante.

On effectue la division de $d$ par $c.$ Le résultat s'écrit
$d=cq+r.$ On a donc $s(d,c;\hat{z}_1)=s(r,c;\hat{z}_1-\hat{q}_1)$
d'après (\ref{stranslation}). On se ramène ensuite au cas où $c_1$
et $r_1$ sont positifs grâce à (\ref{sde-c}) et (\ref{CONJ-D}). La
loi de réciprocité permet alors d'exprimer $s(r,c;.)$ comme une
somme de $s(c,r;.),$ de $s(0,1;.)$ et d'un terme explicite
élémentaire. Cette expression est le résultat final de la boucle.
Il ne reste plus qu'à recommencer en remplaçant le couple $(d,c)$
par le couple $(c,r).$

\smallskip

Puisque $\mc O_F$ est euclidien pour la norme, l'entier $|N_{F/\mb
Q}(r)|$ diminue à chaque division. Par conséquent l'algorithme se
termine, et le dernier reste $\tilde{r}$ est nul. Notons de plus
que l'identité $c\mc O_F+d\mc O_F=\mc O_F$ est préservée à chaque
étape de l'algorithme. Le dernier couple $(\tilde{c},\tilde{r})$
est donc du type $(\tilde{c},\tilde{r})=(\epsilon,0),$ où
$\epsilon$ est une unité de $F.$ Il suffit pour conclure de
transformer la somme de Dedekind $s(0,\epsilon;.)$ en $s(0,1;.)$ à
l'aide de l'identité (\ref{sdcepsi}). La proposition \ref{remeucl}
est ainsi établie dans le cas où $\mc O_F$ est euclidien pour la
norme.

\smallskip
Pour terminer la démonstration, il faut introduire une
généralisation de la notion d'anneau euclidien due à Cooke ([Co]).
Un anneau d'entiers $\mc O_F$ est dit \textit{euclidien en
$k$-étapes pour la norme} si pour tous $c\neq 0$ et $d$ éléments
de $\mc O_F,$ on a besoin de $n\leq k$ divisions successives
$d=cq_1+r_1,$ $ c=r_1q_2+r_2,$ etc... pour obtenir un reste $r_n$
vérifiant $|N_{F/\mb Q}(r_n)|<|N_{F/\mb Q}(c)|.$ En particulier,
les anneaux euclidiens en une étape sont les anneaux euclidiens.

Ceci étant, on a le résultat suivant ([Co, Th.1]) : si $F$ est un
corps de nombres de nombre de classes 1 dont le groupe des unités
est de rang $\geq 1,$ alors $\mc O_F$ est euclidien en $k$-étapes
pour la norme pour un certain entier $k.$

Ce théorème permet de conclure car l'algorithme précédent s'adapte
sans peine au cas d'un anneau d'entiers euclidien en $k$-étapes.
En effet, sous cette hypothèse, il suffit d'effectuer au plus $k$
divisions avant que la norme ne diminue. On utilise la boucle
donnée précédemment à chacune de ces divisions. L'algorithme se
termine donc après avoir effectué au plus $k|N_{F/\mb Q}(c)|$
divisions. On en déduit que l'on peut écrire $s(d,c;\hat{z}_1)$ en
utilisant au plus $k|N_{F/\mb Q}(c)|$ termes du type $s(0,1;.).$
La proposition \ref{remeucl} s'ensuit.
\end{proof}

\smallskip
\noindent Remarque : au cours de l'algorithme précédent, la
fonction $s(0,1;.)$ est évaluée en des $\hat{Z}_1$ qui sont des
images de $\hat{z}_1\in \mc H^{n-1}$ sous l'action de $SL_2(\mc
O_F)$ et des involutions $z_k\mapsto-\bar{z}_k,$
$k\in\{2,\ldots,n\}.$

Nous donnons un exemple qui va rendre la proposition \ref{remeucl}
explicite. On renvoie le lecteur à [Le] pour une bibliographie
très complète sur les corps de nombres euclidiens.

\medskip

\noindent \textbf{Exemple :} Le corps $F=\mb Q(\sqrt{7})$ est
euclidien pour la norme. On va exprimer $s(d,c;z_2)$ en fonction
de $s(0,1;.)$ pour $c_1=3+\sqrt{7}$ et $d_1=-2-\sqrt{7}.$

La division euclidienne peut s'écrire $d=cq+r,$ avec $q=-1$ et
$r=1.$ D'après (\ref{stranslation}), il vient
$s(d,c;z_2)=s(1,c;z_2+1).$

Les hypothèses de la loi de réciprocité étant satisfaites, on
obtient
$$s(1,c;z_2+1)+s(c,1;(\bar{z}_2+1)^{-1})=s(0,1;z_2+1)-\frac1 4 +\kappa_F T(z_2),$$
où $T(z_2)$ est le terme explicite
$$T(z_2)=\left(\frac{1}{3+\sqrt{7}}
+\frac{3+\sqrt{7}}{\left|z_2+1\right|^2} + \frac 1 {(3+\sqrt{7})
\left|(3-\sqrt{7})z_2+4-\sqrt{7}\right|^2} \right)\Ima(z_2).$$
\noindent On note enfin grâce à (\ref{stranslation}) que
$s(c,1;(\bar{z}_2+1)^{-1})=s(0,1;(\bar{z}_2+1)^{-1}-1).$

En mettant bout à bout toutes ces égalités, on trouve en
définitive
\begin{equation}\label{sdcs01exple}
s(d,c;z_2)=s(0,1;z_2+1)-s(0,1;(\bar{z}_2+1)^{-1}-1)-\frac1 4
+\kappa_F T(z_2).
\end{equation}

\bigskip
On se propose maintenant d'étudier plus en détail la fonction
$s(0,1;.)$ qui apparaît dans la proposition précédente. Dans le
cas où $F$ est un corps quadratique réel, la somme de Dedekind
fondamentale $s(0,1;z_2)$ apparaît déjà en filigrane dans le
travail [He] de Hecke. Dans cet article, Hecke étudie une fonction
$\varPsi(z_1,z_2)$ holomorphe sur $ \mc H^2$ qu'il considère comme
un analogue pour $F$ de $\ln \eta.$ Il s'agit en fait du morceau
holomorphe par rapport à $z_1$ \textit{et} $z_2$ de la fonction
$h(z_1,z_2)$ de Asai. En particulier, cette fonction $\varPsi(z)$
est un morceau des fonctions $\Lambda_1$ et $\Lambda_2$ que nous
avons définies précédemment.

Hecke exprime $\varPsi(z)$ à l'aide de fonctions $L$ de Hecke du
corps $F$ via la transformation de Mellin. Il déduit alors de
l'équation fonctionnelle des fonctions $L$ une formule qui relie
$\varPsi(-z^{-1})$ à $\varPsi(z).$

On note $\chi_m$ le Grössencharakter de $F^*$ défini par
$\chi_m(\mu):=\left|\frac{\mu_2}{\mu_1}\right|^{\frac{i\pi m} {\ln
\epsilon}},$ où $m$ est un entier, $\epsilon>1$ est le générateur
du groupe $U_F^+,$ et $\upsilon_1$ le caractère de $F^*$ défini
par $\upsilon_1(\mu):=\sign(\mu_1\mu_2).$ La formule de Hecke fait
apparaître des valeurs spéciales des fonctions
$L(s,\chi_m\upsilon_1)$ de Hecke sur la droite $\Reel(s)=1.$ Ces
résultats se traduisent en termes de somme de Dedekind généralisée
$s(0,1;z_2)$ dans la proposition suivante.

\smallskip

\begin{prop}\label{s01prop}~
Soit $F$ un corps quadratique réel de nombre de classes $1.$ On
note $\lambda_F$ la constante $-\frac{d_F}{2\pi^2R_F^2},$ et on
écrit un élément $z_2$ de $\mc H$ sous la forme $z_2=i \rho
e^{i\frac \pi 2 \theta}$ où $\theta \in \ ]-1,1[$ et $\rho >0.$ La
fonction somme de Dedekind fondamentale $s(0,1;z_2)$ est donnée
par la série :

$\mathrm{a)}$ si toutes les unités de $F$ sont de norme 1 :

$$
s(0,1;i\rho e^{i\frac \pi 2\theta})=\lambda_F\sum_{m=0}^\infty
|L(1+\frac {im\pi}{R_F},\chi_m\upsilon_1)|^2 \frac {\sinh
\left(\frac{m\pi^2}{R_F}\theta \right)} {\sinh\left(\frac {m\pi^2}
{R_F}\right)}\cos\left(2m\pi\frac{\ln \rho}{R_F}\right),$$ où le
terme correspondant à $m=0$ vaut $\frac \theta 2
L(1,\upsilon_1)^2.$

\smallskip

$\mathrm{b)}$ si $F$ a une unité de norme $-1$ :
$$ s(0,1;i\rho
e^{i\frac \pi 2\theta})=\lambda_F\sum_{m=1,\atop m \, \mathrm{
impair}}^\infty |L(1+\frac {im\pi}{2R_F},\chi_m\upsilon_1)|^2
\frac {\sinh \left(\frac{m\pi^2}{2R_F}\theta \right)}
{\sinh\left(\frac {m\pi^2} {2R_F}\right)}\cos\left(m\pi\frac{\ln
\rho}{R_F}\right).$$
\end{prop}

\medskip
\noindent Remarques :

1) La fonction $\varPsi(z)$ se transforme sous l'effet d'une
substitution modulaire quelconque selon une formule donnée dans
[DLT-G, Ex. 2]. On peut obtenir grâce à ces résultats une formule
pour les sommes de Dedekind généralisées $s(d,c,z_2)$ similaire à
celle de la proposition \ref{s01prop}. Elle fait apparaître, comme
dans le cas des sommes de Dedekind classiques, une contribution de
chaque classe $r \! \mod c$ (voir [Ch]).

2) On sait que les valeurs spéciales $L(1+\frac
{im\pi}{R_F},\chi_m\upsilon_1)$ sont non nulles d'après [We, Th.11
p. 288]. Elles interviennent également dans le travail de Arakawa
[Ar].

\smallskip

\noindent \textit{Démonstration de la proposition \ref{s01prop}.}
Comme Hecke, nous définissons pour $\tau_1$ et $\tau_2$ deux
complexes du demi-plan $\Reel(\tau)>0$ la série absolument
convergente
$$ F(\tau_1,\tau_2;\upsilon):= \sum_{\mu\in \mc O_F} \!\!\!\!\begin{array}{l} \prime\\ \
\end{array}\!\!
\left[\sum_{(\nu)|(\mu)}|N_{F/\mb
Q}(\nu)|^{-1}\right]\upsilon(\mu) e^{-\frac{2\pi}
{\sqrt{d_F}}(\tau_1|\mu_1|+\tau_2|\mu_2|)},$$ où $\upsilon$
désigne l'un des deux caractères $\upsilon_0(\mu)=1$ ou
$\upsilon_1(\mu)=\sign(\mu_1\mu_2).$

D'après l'égalité (\ref{defOmegaj}), la fonction $ F$ est reliée à
$\Omega$ par la formule
\begin{equation*} \label{lienlambdaF}
4\Omega(z_1,z_2)= F\left(\frac{z_1} i,\frac {z_2}
i;\upsilon_0\right)-F\left(\frac{z_1} i,\frac{z_2}
i;\upsilon_1\right)+ F\left(\frac{z_1}
i,i\bar{z}_2;\upsilon_0\right)+F\left(\frac{z_1}
i,i\bar{z}_2;\upsilon_1\right).\end{equation*}

\noindent On doit à Hecke ([He, Satz 7]) la formule de
transformation de $ F(\tau_1,\tau_2,\upsilon_0)$ :
\begin{equation*}\label{F0eqfct} F\left(\frac 1 {\tau_1},\frac 1
{\tau_2};\upsilon_0\right)-
F(\tau_1,\tau_2;\upsilon_0)=\zeta_F(2)\frac{\sqrt{d_F}}{\pi^2}\left(\tau_1\tau_2-\frac
1 {\tau_1\tau_2}\right) -\frac{2R_F}{\sqrt{d_F}}\ln
(\tau_1\tau_2).
\end{equation*}

De même, on trouve dans [He, Satz 8] une formule pour la
différence $ F\left(\frac 1 {\tau_1},\frac 1
{\tau_2};\upsilon_1\right)- F(\tau_1,\tau_2;\upsilon_1)$ (noter
qu'il manque un facteur $ \pi^{-1}$ dans le membre de droite de la
dernière égalité p. 402). Cette différence vaut
\begin{equation*}\label{F1eqfct} \frac{\sqrt{d_F}}{\pi
R_F}\left[\frac{2L(1,\upsilon_1)^2}{\pi} \ln (\tau_1\tau_2)+i
\sum_{m\in \mb Z}\!\!\!\begin{array}{l} \prime\\ \ \end{array}
\!\! a_m\left[ \tau_1^{-\frac{2im\pi}{\ln \epsilon}}
+\tau_2^{-\frac{2im\pi}{\ln \epsilon}}\right]\right],
\end{equation*}
où $L(1,\upsilon_1)=0$ par convention s'il existe une unité de
norme $-1,$ et
\begin{equation*}\label{defam}a_m=\left\{ \begin{array}{ll}\dis\frac {|L(1+\frac
{im\pi}{\ln \epsilon},\chi_m\upsilon_1)|^2} {\sinh\left(\frac
{m\pi^2} {\ln \epsilon}\right)} &\textrm{ si } \
\chi_m\upsilon_1(\mu)\
\textrm{ ne dépend que de l'idéal } \ (\mu),\\
0& \textrm{ sinon. } \end{array}\right.\end{equation*}

\noindent Ces deux résultats de Hecke se traduisent en termes de
la fonction $\Omega.$ Si on note $S$ la matrice $
\left(\begin{array}{cc} 0 &-1\\ 1& 0
\end{array}\right),$
on obtient d'après la définition de $\Phi$ :
\begin{equation*} \label{PHISUML}
\Phi\left(S,z_2\right)=i\frac{d_F L(1,\upsilon_1)^2}{4\pi^3R_F^2}
\ln
\left[-\frac{\bar{z}_2}{z_2}\right]-\frac{d_F}{8\pi^2R_F^2}\sum_{m\in\mb
Z}\!\!\!\begin{array}{l} \prime\\ \ \end{array}\!\! a_m
\left[(i\bar{z}_2)^{-\frac{2im\pi}{\ln
\epsilon}}-(-iz_2)^{-\frac{2im\pi}{\ln \epsilon}}\right].
\end{equation*}

\noindent On pose $\lambda_F=-\frac{d_F}{2\pi^2R_F^2}.$ D'après la
définition \ref{definitionsommesgene} et l'égalité précédente, la
fonction $s(0,1;z_2)$ apparaît naturellement comme la partie
réelle de la fonction holomorphe sur $\mc H$ :

$$z_2\mapsto -\frac{\lambda_F} \pi L(1,\upsilon_1)^2\ln (-iz_2)-\frac {\lambda_F}2\sum_{m\in\mb
Z}\!\!\!\begin{array}{l} \prime\\ \ \end{array}\!\! a_m
\left[-iz_2\right]^{\frac{2im\pi}{\ln \epsilon}}.$$

Si on écrit $z_2$ sous la forme $z_2=i\rho e^{i\frac \pi 2\theta}$
avec $\theta$ dans l'intervalle $]-1,1[,$ l'égalité précédente se
réduit à
$$s(0,1;i\rho e^{i\frac \pi 2 \theta})=\lambda_F\frac{\theta} 2
L(1,\upsilon_1)^2+\frac{\lambda_F}{2}\sum_{m\in\mb
Z}\!\!\!\begin{array}{l} \prime\\ \ \end{array}\!\! a_m
\rho^{-\frac{2im\pi}{\ln \epsilon}}\sinh \left(\frac{m\pi^2}{\ln
\epsilon}\theta \right).
$$ Il ne reste plus qu'à constater que $a_{-m}=-a_m$ pour
conclure.\hfill $\Box$

\subsubsection{Opérateurs de Hecke.} Soit $p$ un
nombre premier. On trouve déjà dans le travail de Dedekind [De]
l'identité
\begin{equation*}
s(dp,c)+\sum_{r\, \textrm{mod} \, p} s(d+cr,cp)=(p+1)s(d,c),
\end{equation*}
où $s(.,.)$ désigne la somme de Dedekind classique. Nous nous
proposons dans ce paragraphe de généraliser cette identité dans le
cas d'un corps de nombres $F$ totalement réel de degré $n,$ de
nombre de classes 1.

Pour cela, on désigne par $\mc F$ le $\mb R$-espace vectoriel
formé des fonctions $$f : \mc O_F\times (\mc
O_F\setminus\{0\})\times \mc H^{n}\rightarrow \mb R $$ qui
vérifient la propriété d'invariance par translation :
$$f(d+ cq,c;z+q)=f(d,c;z),$$ où $q$ est un entier de $\mc O_F$ quelconque. On a noté $z+q$ le
$n$-uplet $(z_k+q_k)\in \mc H^n.$

Soit $p$ un entier premier de $\mc O_F$ totalement positif. On
définit \textit{l'opérateur de Hecke} $T_p : \mc F\rightarrow \mc
F $ par la règle
$$(f|T_p)(d,c\,;z):=f(dp,c\,;pz)+\sum_{r\, \mathrm{mod}\, p}
f\left(d+cr,cp\,;\frac{z+r}{p}\right),$$ avec des notations
évidentes pour $pz$ et $\frac{z+r}{p}.$

\noindent Un calcul élémentaire montre que les opérateurs $T_p$
commutent deux à deux.

Soit $\mc F^j$ le sous-espace de $\mc F$ formé des fonctions qui
ne dépendent pas de la variable $z_j\in \mc H.$ La restriction de
$T_p$ à $\mc F^j$ est un endomorphisme de $\mc F^j$ que l'on note
$T_p^j.$

Dans le cas rationnel, la définition de l'opérateur $T_p^1$ est
essentiellement équivalente à celle donnée dans [Ma, §3], où le
lien entre les opérateurs $T_p^1$ et les opérateurs de Hecke
standards est expliqué. En termes d'opérateurs de Hecke $T_p^1$,
l'identité de Dedekind montre que la somme de Dedekind classique
$s(.,.)\in \mc F^1$ est une fonction propre de $T_p^1$ de valeur
propre associée $p+1.$

Dans le cas général, la somme $s_j$ peut être considérée comme un
élément de $\mc F^j$ d'après la proposition \ref{Dsigne}.4. Le
résultat suivant établit la généralisation souhaitée de l'identité
de Dedekind.

\medskip
\begin{prop}\label{sommeoperHecke}
Soit $s_j\in \mc F^j$ la somme de Dedekind généralisée associée au
$j^\textrm{ème}$ plongement réel de $F.$ Pour tout entier $p$
premier de $\mc O_F$ totalement positif on a l'égalité
$$
s_j|T_p^j=\left(N_{F/\mb Q}(p)+1\right)s_j.
$$
Autrement dit, la fonction $s_j$ est une fonction propre pour
$T_p^j$ de valeur propre associée $N_{F/\mb Q}(p)+1.$
\end{prop}

\smallskip\noindent
La preuve de cette proposition résulte immédiatement de la
proposition \ref{propdefsommes} et du lemme suivant :

\begin{lem}\label{propOmegap}
Soit $p\in \mc O_F$ un premier totalement positif, et $z\in \mc
H^n.$ Alors $$ \Omega_j(pz)+ \sum_{r\, \mathrm{mod}\, p}
\Omega_j\left(\frac{r+z}{p}\right)=\left(N_{F/\mb
Q}(p)+1\right)\Omega_j(z).
$$
\end{lem}

\smallskip
\noindent \textit{Démonstration du lemme \ref{propOmegap} :} La
fonction $\Omega_j$ étant donnée par (\ref{defOmegaj}), on
considère la somme $\sum_{r\, \textrm{mod}\, p}
\Omega_j\left(\frac{r+z}{p}\right)$. Cette somme est bien définie
car $\Omega_j(z)$ est invariante quand on translate $z$ par un
entier de $\mc O_F.$ Elle vaut
$$\sum_{\nu \in \mc{O}_F/U_F^+} \!\!\!\prim\!\!\!\!
\frac{[U_F:U_F^+]^{-1}}{|N_{F/\mb Q}(\nu)|} \sum_{ \mu \in
\mc{O}_F, \atop \frac{\mu_j \nu_j}{ \mf{\delta}_j}>0}\!\prim
\!\!\!\!\!\!\! \left(e^{2i\pi\sum_{k=1}^n
\left(\frac{\mu_k\nu_k}{p_k\delta_k}
x_k+\left|\frac{\mu_k\nu_k}{p_k\delta_k}\right|iy_k\right)}\right)\sum_{r\,
\textrm{mod}\, p}e^{2i\pi\, \Tra_{F/\mb Q} \left(\frac{\mu\nu
r}{p\mf{\delta}}\right)}.
$$

On remarque alors que $$\sum_{r\, \textrm{mod}\, p}e^{2i\pi\,
\Tra_{F/\mb Q} \left(\frac{\mu\nu
r}{p\mf{\delta}}\right)}=\left\{\begin{array}{lll}N_{F/\mb Q}(p)
&\ \ \textrm{si}\ p \ \textrm{divise}\ \mu \nu, \\ 0 & \ \
\textrm{sinon.}
\end{array}\right.
$$
Puisque $p$ est premier, on obtient une décomposition de la forme
suivante :
$$\sum_{r\, \textrm{mod}\, p}
\Omega_j\left(\frac{r+z}{p}\right)=\sum_{\dis \nu,\mu \atop \dis
p\,|\nu}\ +\ \sum_{\dis \nu,\mu\atop \dis p\,|\mu}\ -\! \!
\sum_{\dis\nu,\mu\atop\dis p\,|\nu \ \mathrm{ et }\,
p\,|\mu}.$$\noindent La première de ces trois sommes est
$\Omega_j(z),$ la deuxième est $N_{F/\mb Q}(p)\Omega_j(z),$ et la
troisième $\Omega_j(pz).$ Le résultat s'ensuit.\hfill $\square$

\section{Valeur spéciale de fonctions $L$ et invariants de classes de
$SL_2(\mc O_F).$}\label{partieL} \noindent Dans la partie
précédente, nous avons défini et étudié les sommes de Dedekind
généralisées $s_j(d,c;.): \mc H^{n-1}\rightarrow \mb R$. Elles
sont construites, à un terme élémentaire près, à partir de la
généralisation $\Phi_j\left(\left(
\begin{array}{cc}a&b\\c&d\end{array} \right),\,.\right): \mc
H^{n-1}\rightarrow \mb R$ de la fonction $\Phi_R$ de Rademacher.

Ceci étant, nous nous intéressons dans le reste de cet article à
la valeur spéciale $\Phi_j(A,\omega_c)$ de la fonction
$\Phi_j(A,.)$ en un point algébrique $\omega_c\in \mc H^{n-1}$
associé à un certain type de matrices $A$ de $\Gamma.$

\noindent On rappelle qu'une matrice $M$ de $SL_2(\mb R)$ est dite
elliptique si $|\tra(M)|<2,$ hyperbolique si $|\tra(M)|>2.$
Autrement dit, une matrice elliptique a un unique point fixe dans
$\mc H,$ alors qu'une matrice hyperbolique a deux points fixes à
l'infini. Enfin, une matrice $(A_k)\in SL_2(\mb R)^n$ est dite
elliptique si chacune de ses composantes est elliptique.

Considérons une matrice $A$ de $\Gamma$ comme un élément $(A_k)$
de $SL_2(\mb R)^n.$ Soit $j\in\{1,\ldots,n\}$ un entier fixé. On
suppose que le $(n-1)$-uplet de matrices $(A_k)_{k\neq j}$ a un
point fixe $\omega_c\in \mc H^{n-1}.$ On normalise la valeur
spéciale $\Phi_j(A,\omega_c)$ de la façon suivante.

\begin{defi} \label{defPsi}Soit $A\in \Gamma$ et $j\in\{1,\ldots,n\}$ comme ci-dessus. On définit le réel $\Psi_j(A)$ en
posant \begin{equation*}
\Psi_j(A):=2^nR_F\Phi_j(A,\omega_c)-2^{n-2}R_F\sign(c_j\tra(A_j)),\end{equation*}
où $R_F$ désigne le régulateur du corps $F.$ Considérons en
particulier le cas d'une matrice $A\in\Gamma$ ayant une seule
composante $A_j$ hyperbolique. Ses autres composantes sont alors
elliptiques. Par suite, le réel $\Psi_j(A)$ est bien défini et
sera noté $\Psi(A)$ sans ambiguïté.
\end{defi}

\smallskip

On verra que $\Psi_j(A)$ est un invariant de classes de $\Gamma$
de nature différente selon que la composante $A_j$ de $A$ est
hyperbolique ou elliptique. Dans le cas où $A_j$ est hyperbolique,
nous allons tout d'abord définir une certaine fonction holomorphe
$L_A(s)$ dont la valeur spéciale en $s=0$ sera reliée à $\Psi(A)$
dans le théorème \ref{thL0PSI}. Notre exposé trouve sa source dans
les résultats de Hara ([Ha]). Notre présentation s'inspire celle
de Atiyah ([At]) pour le groupe $SL_2(\mb Z).$

\subsection{Valeur spéciale de fonctions $L$ en $s=0$.}
\label{fctionL} Il est nécessaire d'établir quelques résultats
préliminaires sur les matrices de $\Gamma$ ayant une seule
composante hyperbolique avant de définir la fonction $L$ qui leur
est associée.
\subsubsection{Généralités sur les matrices
quasi-elliptiques.}\label{partiequasiell}
\begin{defi}\label{defquasell}
Une matrice $A$ de $\Gamma$ sera dite \textit{quasi-elliptique}
si, considérée comme une matrice de $SL_2(\mb R)^n,$ elle a une
seule composante $A_j$ hyperbolique. Les autres composantes sont
alors toutes elliptiques, et on note $\omega_c\in \mc H^{n-1}$ le
point fixe de $(A_k)_{k\neq j}.$

On dira qu'une extension quadratique $K$ de $F$ est
\textit{quasi-totalement complexe} (resp. \textit{quasi-totalement
réelle}) si elle a $n-1$ places complexes (resp. 1 place
complexe).
\end{defi}

\smallskip

Si $A$ est une matrice quasi-elliptique, ses valeurs propres
engendrent une extension quasi-totalement complexe $K=K_A$ de $F.$
Le lemme suivant résulte immédiatement de la définition
\ref{defquasell} et du th\'eor\`eme des unit\'es de Dirichlet. Sa
démonstration est laissée au lecteur.

\begin{lem} \label{lemquasiell}
Soit $A=\left(
\begin{array}{cc}a&b\\c&d\end{array} \right)$ une matrice de
$\Gamma$ quasi-elliptique et $K=K_A$ l'extension quadratique de
$F$ associée. Alors :
\begin{enumerate}
\item L'entrée $c$ de $A$ est non nulle. \item Le groupe des
unités $U_K$ de $K$ est de rang $n.$ Ainsi on a $[U_K:U_F]=1.$
\item Les valeurs
propres de $A$ s'écrivent $ {1\over 2}
\left[\tra(A)\pm\sqrt{\tra(A)^2-4}\right].$ Ce sont des unités
relatives de l'extension $K/F.$ \item Soit
$\Gamma_{\omega_c}:=\{B\in \Gamma / (B_k)_{k\neq
j}\,\omega_c=\omega_c\}$ le sous-groupe de $\Gamma$ constitu\'e
des matrices quasi-elliptiques ayant m\^eme point fixe que $A.$ Il
est isomorphe $\{\pm 1\}\times \mb Z.$ En outre, si $P$ est une
matrice de $\Gamma,$ alors  $P^{-1}AP$ est quasi-elliptique et on
a $P\Gamma_{\omega_c}P^{-1}=\Gamma_{P\omega_c}$ et
$K_A=K_{P^{-1}AP}.$
\end{enumerate}

\end{lem}
\smallskip

\noindent \textbf{ Notations.}

Soit $g$ le générateur du groupe de Galois de l'extension $K/F.$
On note $\{\omega,\omega^g\}\in K ,$ et on appelle \textit{points
fixes} de la matrice quasi-elliptique $A$ les deux générateurs de
l'extension $K/F$ solutions de l'équation
\begin{equation*}\label{defomega}(cX+d)X=aX+b.\end{equation*}

Soit $\omega\in K$ un des deux points fixes de $A,$ choisi
arbitrairement. Quitte à renuméroter les plongements réels de $F,$
on peut supposer que $A_1$ est la composante hyperbolique de $A.$
On ordonne les deux plongements réels $r_1,r_2$ de $K$ au-dessus
de $\iota_1$ de façon à ce que $\omega_{r_2}<\omega_{r_1};$ on
note encore $\iota_k$ le plongement complexe de $K$ qui prolonge
le plongement réel $\iota_k$ de $F$ de façon à ce que le point
fixe $\omega_c\in \mc H^{n-1}$ de $(A_k)_{k=2}^n$ soit le
$(n-1)$-uplet $(\omega^{\iota_2},\ldots,\omega^{\iota_n}).$
L'image $\beta^{\iota_k}$ d'un élément $\beta\in K$ par $\iota_k$
sera notée $\beta_k.$

\noindent Les deux valeurs propres de $A$ s'écrivent
$\epsilon=c\omega+d$ et $\epsilon^g=\epsilon^{-1}=c\omega^g+d.$

Nous introduisons un ordre $\mc O_A$ de $K,$ un groupe d'unités
$U_A\subset \mc O^*_A$ et un module $\mc M_\omega\subset K $ sur
l'ordre $\mc O_A$ qui vont nous permettre de définir la fonction
$L_A(s).$

On pose $\mc O_A=\mc O_F+\epsilon\mc O_F\subset \mc O_K.$ C'est
bien un ordre de $K$ contenant $\epsilon$ et $\epsilon^{-1}$ en
vertu de l'égalité $\epsilon+\epsilon^{-1}=\tra(A)$ et du lemme
\ref{lemquasiell}.3. En outre, le lemme \ref{lemquasiell}.4 montre
que $\epsilon$ est sans torsion dans $U_K/U_F,$ ce qui permet de
définir le produit direct $U_A:=U_F\times <\epsilon>.$ Ce
sous-groupe d'indice fini du groupe des unités de $K$ est contenu
dans l'ordre $\mc O_A.$

On définit enfin le $\mc{O}_F$-module libre de rang 2 $\mc
M_\omega :=\mc{O}_F+\omega\mc{O}_F\subset K.$ C'est un module sur
l'ordre $\mc{O}_ A$ car on a $\epsilon\omega=a\omega+b.$ Il ne
dépend pas seulement de la matrice $A$ mais aussi du choix de
$\omega.$

\subsubsection{La fonction $L_A(s)$ associée à une matrice
quasi-elliptique.}

Les notations sont celles du paragraphe précédent.

\begin{defi}\label{textdefLA}
Soit $A$ une matrice quasi-elliptique dont la composante $A_1$ est
hyperbolique. On lui associe la fonction holomorphe $L_A(s)$
définie pour $\Reel(s)>1$ par la formule
\begin{equation*}\label{defLA}L_A(s):= \sign(c_{1}\tra(A_{1}))
\sum_{\beta \in \mc M_\omega /U_A}\!\!\!\!\prim\!\!\!
\frac{\sign(\beta_{r_1}\beta_{r_2})}{|N_{K/\mb
Q}(\beta)|^s}.\end{equation*}
\end{defi}

Remarques :
\begin{itemize}
\item Cette fonction ne dépend effectivement que de $A$ et
pas du choix du point fixe $\omega$ de $A.$ En effet, quand on
change $\mc M_\omega$ en $\mc M_{\omega^g},$ on ne change pas la
fonction $L$ car $\beta\mapsto \beta_{r_1}\beta_{r_2}$ et
$\beta\mapsto N_{K/\mb Q}(\beta)$ sont deux fonctions
$g$-invariantes. \item Nous démontrerons dans le corollaire
\ref{LAEth} que $L_A(s)$ admet un prolongement sur $\mb C$ en une
fonction entière qui vérifie une équation fonctionnelle.

\item Dans le cas particulier où $\mc M_\omega$ est un $\mc
O_K$-module, la fonction est $L_A$ est simplement la fonction $L$
de Hecke partielle associée à l'idéal $\mc M_\omega$ de $\mc O_K$
et au caractère $\varphi : \beta \mapsto
\sign(\beta_{r_1}\beta_{r_2})$ de $K^*.$

\end{itemize}

\smallskip La proposition suivante permet de comparer les fonctions $L_{P^{-1}AP}$ et $L_{A^k}$ à la fonction $L_A.$
En particulier, elle justifie la présence du facteur
$\sign(c_1\tra(A_1))$ dans la définition de $L_A$ pour assurer
l'identité $L_{A^{-1}}=-L_A.$

\smallskip

\begin{prop}\label{propLAPAP} Soit $A$ une matrice
quasi-elliptique de $\Gamma.$
\begin{enumerate}
\item Soit $P=\left(\begin{array}{cc} a'' &b''\\
c''& d''
\end{array}\right)\in \Gamma.$ Alors
$$L_{PAP^{-1}}(s)=|N_{K/\mb Q}(c''\omega+d'')|^s L_A(s).$$
\item Soit $k$ un entier relatif non nul. Alors $$L_{\pm
A^k}(s)=kL_A(s).$$
\end{enumerate}
\end{prop}

\begin{proof}

Pour établir la première assertion, nous allons étudier
successivement le comportement de $c,$ $ U_A$ et $\mc M_\omega$
quand on change $A$ en $PAP^{-1}.$

\smallskip

\begin{itemize}

\item[i)]
En notant $A'= PAP^{-1}=\left(\begin{array}{cc} a' &b'\\ c'& d'
\end{array}\right),$ un calcul direct établit les égalités suivantes
entre les valeurs propres de $A$ et celles de $A'$ :
$$ c\omega+d=c'(P\omega)+d'\ \ \ \ \ \mathrm{ et
}\ \ \ \ \ c\omega^g+d=c'(P\omega^g)+d'.$$

D'où par soustraction des deux égalités précédentes
$$c(\omega-\omega^g)=c'\left(\frac{a''\omega+b''}{c''\omega+d''}-\frac{a''\omega^g+b''}{c''\omega^g+d''}\right).$$
Après avoir chassé les dénominateurs, on obtient enfin
\begin{equation}c'=c
N_{K/F}(c''\omega+d'').\label{signetrA1}\end{equation}

\item[ii)] Les matrices $PAP^{-1}$ et $A$ ont mêmes valeurs
propres $\{\epsilon,\epsilon^{-1}\},$ donc
\begin{equation}\label{UPAP} U_{PAP^{-1}}=U_A.\end{equation}
\item[iii)]

Par définition, $\{\omega,1\}$ est une base de $\mc M_\omega$
relativement à $\mc O_F.$ Puisque $P$ appartient à $\Gamma,$
$\{a''\omega+b'',c''\omega+d''\}$ est une autre base. Ainsi
$$\mc M_\omega=(a''\omega+b'')\mc O_F+(c''\omega+d'')\mc
O_F=(c''\omega+d'')(\mc O_F+(P\omega)\mc O_F),$$ qui s'écrit aussi
\begin{equation}\label{MPAP}\mc M_{P\omega}=\frac 1
{c''\omega+d''} \mc M_\omega.\end{equation}\end{itemize}

\smallskip\noindent
Les points fixes de $PAP^{-1}$ sont $P\omega$ et $P\omega^g.$ En
utilisant (\ref{signetrA1}) et (\ref{UPAP}), il vient
\begin{equation*}L_{PAP^{-1}}(s)=\sign(c_1\tra(A_1)N_{K/F}(c''\omega+d'')_1)\sum_{\beta\in
\mc M_{P\omega}/U_A}\!\!\!\!\!\prim \!
\frac{\sign(N_{K/F}(\beta)_1)}{|N_{K/\mb
Q}(\beta)|^s}.\end{equation*} \noindent On conclut à l'égalité
souhaitée grâce à (\ref{MPAP}).

Démontrons maintenant la deuxième assertion.

\smallskip
L'identité $L_{\pm A}=L_A$ résulte immédiatement de la définition
de $L_A(s).$

Le cas où $k=-1$ est également simple. Les matrices $A$ et
$A^{-1}$ ont même trace, mêmes valeurs propres
$\{\epsilon,\epsilon^{-1}\}$ et mêmes points fixes
$\{\omega,\omega^g\}.$ Leurs entrées $c$ étant de signe opposé, il
s'ensuit que $L_{A^{-1}}(s)=-L_A(s).$

Il suffit maintenant de montrer que $L_{A^k}=kL_A$ si $k$ est un
entier positif. Les valeurs
propres de $A^k=\left(\begin{array}{cc} a^{(k)} &b^{(k)}\\
c^{(k)}& d^{(k)}
\end{array}\right)$ sont $\{\epsilon^k,\epsilon^{-k}\},$ et ses points fixes sont $\{\omega,\omega^g\}.$ Ainsi le groupe
$U_{A^k}$ est un sous-groupe d'indice $k$ de $U_A,$ et par suite
$kL_A$ et $L_{A^k}$ coïncident au signe près. Il reste à régler le
problème du signe de $c_1^k\tra(A_1^k).$

\smallskip
On a convenu que $\omega_{r_2}<\omega_{r_1}.$ Il en résulte que
pour tout $k\geq 1,$ on a
\begin{equation}\label{signcepsi}\epsilon_{r_1}^k=\frac
12\left[\tra(A_1^k)+\sign(c_1^{(k)})\sqrt{\tra(A_1^k)^2-4}\right].\end{equation}
\noindent D'après l'égalité précédente, si $c_1^{(k)}\tra(A_1^k)$
est positif alors $|\epsilon_{r_1}^k|>|\tra(A_1^k)|/ 2.$ On en
déduit que $|\epsilon_{r_1}^k|>1$ car $A_1^k$ est hyperbolique. La
réciproque s'obtient en changeant $\epsilon$ en $\epsilon^g.$ On a
donc la série d'assertions équivalentes suivantes :
\begin{equation*}\label{epsisign}\sign(c_1^{(k)}\tra(A_1^k))=+1 \Longleftrightarrow
|\epsilon_{r_1}^k|>1 \Longleftrightarrow |\epsilon_{r_1}|>1
\Longleftrightarrow \sign(c_1\tra(A_1))=+1.\end{equation*}

Ainsi $c_1\tra(A_1)$ et $c_1^{(k)}\tra(A_1^k)$ ont même signe et
on a en définitive
\begin{equation*}L_{A^k}(s)=kL_A(s).\end{equation*}
\end{proof}

\subsubsection{Lien entre la valeur spéciale $L_A^{(n-1)}(0)$ et $\Psi(A).$} Soit $A$ une matrice quasi-elliptique de
$\Gamma$ dont la composante hyperbolique est $A_1.$ Soit
$\omega_{r_2}<\omega_{r_1}$ les deux points fixes réels de $A_1,$
et $\omega_c\in \mc H^{n-1}$ le point fixe des autres composantes
$(A_k)_{k=2}^n.$

Nous avons montré dans le paragraphe précédent comment associer à
la matrice $A$ une fonction holomorphe $L_A(s)$ définie pour
$\Reel(s)>1.$ De plus, la définition \ref{defPsi} permet
d'associer à la matrice $A$ le réel $\Psi(A)$ en renormalisant la
valeur spéciale $\Phi_1(A,\omega_c).$

Le but de ce paragraphe est de prouver d'une part que $L_A(s)$
admet un prolongement holomorphe sur $\mb C$ qui a un zéro d'ordre
au moins $n-1$ en $s=0,$ et d'autre part que la valeur
$L_A^{(n-1)}(0)$ coïncide avec $(n-1)!\,\Psi(A).$ Ces deux
résultats se déduisent d'une formule qui exprime $L_A(s)$ en
fonction de la période d'une forme différentielle $A_1$-invariante
construite à partir la série d'Eisenstein non-holomorphe
$E_F(z,s)$ introduite par Asai.

\smallskip

Plus précisément, en nous inspirant des travaux de Hara [Ha], C.
Meyer [Me] et Siegel [Si] qui utilisent une méthode due
originellement à Hecke, nous considérons la forme différentielle
$$\frac{\partial }{\partial z_1}E_F(z_1,\omega_c,s)dz_1.$$
Il résulte immédiatement de la proposition \ref{propEasai}.2 que
cette forme est $A_1$-invariante. Le demi-cercle $\mc C$ de
diamètre $[\omega_{r_2},\omega_{r_1}]$ est l'unique géodésique de
$\mc H$ invariante par la matrice $A_1.$

Le théorème suivant relie la période de cette forme différentielle
à la fonction $L_A(s).$


\medskip

\begin{theo}\label{EsLA} Soit $\tau$ un point du demi-cercle $\mc C.$ Soit $s$ un complexe de partie réelle $>1.$
Alors on a l'égalité
\begin{equation}
\label{ELA}\int_{\tau}^{A_1\tau} \frac{\partial }{\partial
z_1}E_F(z_1,\omega_c,s)dz_1 =\frac{\Gamma\left( \frac{s+1}
2\right)^2}{\Gamma(s)} \frac{\vol(\mc M_\omega)^s}{2i(d_F)^s}
L_A(s),
\end{equation}
où l'intégrale porte sur l'arc de $\mc C$ joignant $\tau$ à
$A_1\tau.$ Le terme $ \vol(\mc M_\omega)$ désigne le volume du
$\mc O_A$-module $\mc M_\omega$ considéré comme un réseau de $\mb
R^2\times \mb C^{n-1}.$
\end{theo}

\begin{proof}
\smallskip

Soit $\tau$ un point du demi-cercle, et $\tau^*=A_1\tau\in \mc C$
son image par $A_1.$ On pose $$J_A(s):= \int_{\tau}^{\tau^*}
\frac{\partial }{\partial z_1}E_F(z_1,\omega_c,s)dz_1,$$
l'intégrale portant sur l'arc $\mc C_\tau$ joignant $\tau$ à
$\tau^*.$

Nous allons construire une paramétrisation de l'arc $\mc C_\tau.$
Pour cela, définissons l'application $f: z_1\mapsto i\frac
{z_1-\omega_{r_2}}{z_1-\omega_{r_1}}.$ On note que $f$ réalise une
bijection du demi-cercle $\mc C$ sur $\mb R^*_+.$ En outre, il
suffit d'utiliser les égalités $A_1\omega_{r_1}=\omega_{r_1}$ et
$A_1\omega_{r_2}=\omega_{r_2}$ pour vérifier que
$f(\tau^*)=\epsilon_{r_1}^2f(\tau).$

On peut supposer sans perte de généralité que $f(\tau)=1$ et
$|\epsilon_{r_1}|>1.$ En effet, $J_A(s)$ est la période d'une
forme différentielle $A_1$-invariante et ne dépend donc pas du
choix du point base $\tau$ de $\mc C.$ On voit également que,
quitte à changer $A$ en $A^{-1},$ ce qui a pour effet de changer
$J_A(s)$ et $L_A(s)$ en leur opposé, il suffit de démontrer le
théorème dans le cas où $|\epsilon_{r_1}|>1$.

On supposera dans la suite que ces deux hypothèses
simplificatrices sont satisfaites. L'application $g=f^{-1}$
définit alors une paramétrisation de l'arc $\mc C_\tau$ :
$$ g : [1,\epsilon_{r_1}^2] \rightarrow \mc C_\tau.$$

Cela conduit à effectuer dans $J_A(s)$ le changement de variable
$$z_1=g(t)=\frac{t\omega_{r_1}-i\omega_{r_2}}{t-i},\ \ \ \ \ \ \textrm{ d'où }\ \ \
\ \ \ \
dz_1=\frac{(z_1-\omega_{r_1})(z_1-\omega_{r_2})}{\omega_{r_2}-\omega_{r_1}}\frac{dt}t.$$

On calcule explicitement la dérivée partielle $ \frac{\partial
}{\partial z_1}=\frac{\partial }{\partial x_1} -i\frac{\partial
}{\partial y_1}$ de $E_F(z,s)$ à partir de la définition
\ref{defEisnonholo}. On vérifie que $\frac{\partial }{\partial
z_1}E_F(z_1,\omega_c,s)$ est donné pour $\Reel(s)>1$ par la série
\begin{equation}
\label{defformE} \frac{\partial }{\partial z_1}E_F(z_1,\omega_c,s)
=\frac s {2i}\!\!\! \sum_{\ \ \ (\mu,\nu)\in \mc O_F^2/U_F}
\!\!\!\!\!\!\!\!\!\prim\! \frac{ y_1^{s-1}(\mu_1\bar{z}_1+\nu_1)^2
}{|\mu_1z_1+\nu_1|^{2s+2}} \prod_{k=2}^n \frac{(\Ima \
\omega_k)^s}{|\mu_k\omega_k+\nu_k|^{2s}}.\end{equation}

Après un calcul élémentaire, on arrive à
\begin{equation}\label{morcint}\frac{
y_1^{s-1}(\mu_1\bar{z}_1+\nu_1)^2}{|\mu_1z_1+\nu_1|^{2s+2}}dz_1=i
(\omega_{r_1}-\omega_{r_2})^s
\frac{t^s(-i\beta_{r_1}t+\beta_{r_2})^2}{(\beta_{r_1}^2t^2+\beta_{r_2}^2)^{s+1}}\frac{dt}
t, \end{equation} où l'on a posé $\beta=\mu\omega+\nu.$

Lorsque le couple $(\mu,\nu)$ parcourt l'ensemble des classes non
nulles du quotient $\mc O_F^2\mod U_F,$ on observe que $\beta$
parcourt l'ensemble des classes non nulles du quotient $\mc
M_\omega \mod U_F.$ Nous déduisons de (\ref{defformE}),
(\ref{morcint}) et de la remarque précédente que
\begin{equation}\label{Egpresq}J_A(s)=\frac{sV^s}{2} \sum_{\beta\in \mc M_\omega/U_F}\!\!\!\!\prim\!\!\!
\left(\prod_{k=2}^n|\beta_k|^{-2s}\right)\int_{1}^{\epsilon_{r_1}^2}
\frac{t^s(-i\beta_{r_1}t+\beta_{r_2})^2}{(\beta_{r_1}^2t^2+\beta_{r_2}^2)^{s+1}}\frac{dt}t,
\end{equation}
où $V=(\omega_{r_1}-\omega_{r_2})\prod_{k=2}^n \Ima(\omega_k).$ Le
réel $Vd_F$ peut s'interpréter comme le volume $\vol(\mc
M_\omega)$ du module $\mc M_\omega$ vu comme un réseau de $\mb
R^{2}\times \mb C ^{n-1}$ d'après [Sa, Prop. 4.2.1 et 4.2.2].

Utilisons maintenant une idée due à Hecke. Elle consiste d'abord à
effectuer dans l'intégrale du membre de droite de (\ref{Egpresq})
le changement de variable
$u:=\left|\frac{\beta_{r_1}}{\beta_{r_2}}\right|t.$ Puisque
$|N_{K/\mb Q}(\beta)|=|\beta_{r_1}\beta_{r_2}|\prod_{k=2}^n
|\beta_k|^2,$ l'égalité (\ref{Egpresq}) devient
\begin{equation*}\label{ETAPg}
J_A(s)= \frac{s\vol(\mc M_\omega)^s}{2(d_F) ^s}\sum_{\beta \in \mc
M_{\omega}/ U_F}\!\!\!\!\prim\!\!\!|N_{K/\mb
Q}(\beta)|^{-s}\int_{\left|\frac{\beta_{r_1}}{\beta_{r_2}}\right|}^{\left|\frac{\beta_{r_1}}{\beta_{r_2}}\right|
\epsilon_{r_1}^{2}} \frac {u^s(1-i\varphi(\beta)u)^2
}{(u^2+1)^{s+1}} \frac{du} u,\end{equation*} où on a posé
$\varphi(\beta):=\sign(\beta_{r_1}\beta_{r_2}).$

On observe ensuite que l'on obtient un système de représentants
des classes de $\mc M_\omega\setminus\{0\}\mod U_F$ en considérant
la famille $\{\epsilon^k\beta\}$ où $k\in \mb Z$ et $\beta\in \mc
M_\omega\setminus\{0\}\mod U_A.$ Cela nous permet d'écrire
l'égalité
$$J_A(s)=\frac{s\vol(\mc M_\omega)^s}{2(d_F)^s} \sum_{\beta \in \mc
M_{\omega}/ U_A}\!\!\!\!\prim\!\!\!|N_{K/\mb Q}(\beta)|^{-s}
\sum_{k \in \mb
Z}\int_{\left|\frac{\beta_{r_1}}{\beta_{r_2}}\right|\epsilon_{r_1}^{2k}}^{\left|\frac{\beta_{r_1}}{\beta_{r_2}}\right|
\epsilon_{r_1}^{2k+2}} \frac {u^s(1-i\varphi(\beta)u)^2
}{(u^2+1)^{s+1}} \frac{du} u.$$

La famille de segments $\left\{
\left[\left|\frac{\beta_{r_1}}{\beta_{r_2}}\right|\epsilon_{r_1}^{2k},\left|\frac{\beta_{r_1}}{\beta_{r_2}}\right|
\epsilon_{r_1}^{2k+2} \right[,k \in \mb Z \right\}$ forme une
partition de $\mb R_+^*.$ Par conséquent
\begin{equation} \label{avantreportI}J_A(s)=\frac{s\vol(\mc M_\omega)^s}{2(d_F)^s}\sum_{\beta
\in \mc M_{\omega}/ U_A}\!\!\!\!\prim\!\!\!|N_{K/\mb
Q}(\beta)|^{-s} \int_0^\infty \frac {u^s(1-i\varphi(\beta)u)^2
}{(u^2+1)^{s+1}} \frac{du} u.\end{equation}

Pour conclure, calculons la valeur de l'intégrale
$$I:= \int_0^\infty \frac {u^s(1-i\varphi(\beta)u)^2 }{(u^2+1)^{s+1}}
\frac{du} u.$$ En effectuant le changement de variable
$\tilde{u}:=u^{-1},$ on trouve \begin{align*}I=&\int_0^\infty
\frac {u^s }{(u^2+1)^{s+1}}(-i\varphi(\beta)+u)^2 \frac{du}
u\\=&-I-4i \varphi(\beta)\int_0^\infty \frac{u^s}{(u^2+1)^{s+1}}
\frac{du}{u},\end{align*} d'où
\begin{equation*}\label{valeurI}I=-i\frac{\Gamma\left( \frac{s+1}
2\right)^2}{\Gamma(s+1)}\sign(\beta_{r_1}\beta_{r_2}).\end{equation*}

On reporte la valeur de l'intégrale $I$ dans le membre de droite
de (\ref{avantreportI}). Puisque $|\epsilon_{r_1}|>1,$ nous savons
que le facteur $\sign(c_1\tra(A_1))$ qui apparaît dans la
définition de $L_A(s)$ est égal à $1.$ Le théorème s'ensuit.
\end{proof}


\medskip
D'après le théorème précédent, il suffit d'étudier les propriétés
de la forme $A_1$-invariante $\frac{\partial}{\partial
z_1}E_F(z_1,\omega_c,s)dz_1$ pour obtenir le prolongement
holomorphe de $L_A(s).$ Ces propriétés sont rassemblées dans la
proposition suivante qui est une conséquence immédiate du théorème
\ref{propEasai}.
\begin{prop}\label{propEpart} Soit $z\in \mc H^n$ fixé.
\begin{enumerate}
\item Définie pour $\Reel(s)>1$ par une série, la fonction
$s\mapsto \frac{\partial}{\partial z_1}E_F(z,s)$ se prolonge en
une fonction holomorphe sur tout le plan complexe. \item On a
l'équation fonctionnelle $$G_F(2s)\frac{\partial}{\partial
z_1}E_F(z,s)=G_F(2-2s).$$ \item
\smallskip Au
voisinage de $s=0,$ on a le développement en série de Laurent
\begin{equation*} \label{Epartvoisinage0}\frac{\partial}{\partial z_1}E_F(z,s)=-2^{n-2}R_Fs^n\left[
\frac{1}{z_j-\bar{z}_1}-\frac{\partial}{\partial z_1}
h(z)\right]+O(s^{n+1}).
\end{equation*}
\end{enumerate}
\end{prop}

\begin{cor}\label{LAEth}

Pour toute matrice $A$ quasi-elliptique de $ \Gamma$, la fonction
$ L_A$ définie pour $\Reel(s)>1$ se prolonge en une fonction
holomorphe sur tout le plan complexe. De plus la fonction
$$\tilde{L}_A(s):=\left(\vol(\mc M_\omega)\pi^{-n}\right)^s\Gamma(s)^{n-1}\Gamma\left(\frac{s+1}
2\right)^2 L_ A(s)$$ est une fonction holomorphe sur $\mb C,$
invariante quand on change $s$ en $1-s.$
\end{cor}

\noindent \textit{Démonstration du corollaire \ref{LAEth}.} On
déduit du théorème \ref{EsLA} et de la proposition
\ref{propEpart}.1 que $L_A(s)$ se prolonge en une fonction
méromorphe sur $\mb C.$ Ses pôles, s'il y en a, ne peuvent
provenir que des pôles de $\Gamma(s).$ Ils ne peuvent donc se
trouver qu'en $s=-m,$ $m\in \mb N.$ En outre, la proposition
\ref{propEpart}.3 montre que $L_A (s)/s^{n-1}$ est bornée au
voisinage de 0. Par conséquent, la fonction $\tilde{L}_A$ est
holomorphe au moins sur le demi-plan $\Reel(s)>-1.$

\smallskip

On déduit du théorème \ref{EsLA} et de la proposition
\ref{propEpart}.2 l'égalité
$$\tilde{L}_A(s)=\tilde{L}_A(1-s).$$

\noindent Par suite, la fonction $\tilde{L}_A$ est holomorphe sur
les demi-plans $\Reel(s)>-1$ et $\Reel(s)<2,$ donc sur $\mb C.$ On
en conclut que $L_A$ est holomorphe sur $\mb C.$\hfill $\Box$

\medskip

Etudions maintenant de manière plus précise le comportement de
$L_A(s)$ au voisinage de $s=0$ au moyen de la formule limite de
Kronecker généralisée.

\medskip

\begin{theo}\label{thL0PSI}
Pour toute matrice $A\in \Gamma$ quasi-elliptique, la fonction
$L_A$ a un zéro d'ordre supérieur ou égal à $n-1$ en $s=0.$ En
outre on a l'identité
\begin{equation*} L_A^{(n-1)}(0)=(n-1)!\, \Psi(A).\end{equation*}
\end{theo}

\begin{proof}

En utilisant les développements en série de Laurent au voisinage
de $s=0$ donnés par la proposition \ref{propEpart}.3 d'une part,
et $\Gamma(s)=s^{-1} +O(1)$ d'autre part, on déduit du théorème
\ref{EsLA} l'identité
$$
L_A(s)= \frac{2^{n}R_F}{2i\pi }s^{n-1} \int_{\tau}^{\tau^*}\left(
\frac{1}{z_1-\bar{z}_1}-\frac{\partial}{\partial z_1}
h(z_1,\omega_c)\right)dz_1 +O(s^n).$$

On voit donc que la fonction $L_A(s)$ possède donc un zéro d'ordre
$\geq n-1$ en $s=0,$ et que sa dérivée d'ordre $n-1$ en $s=0$
s'écrit
$$L_A^{(n-1)}(0)= \frac{(n-1)!}{2i\pi
}2^{n}R_F \int_{\tau}^{\tau^*}\left(
\frac{1}{z_1-\bar{z}_1}-\frac{\partial}{\partial z_1}
h(z_1,\omega_c)\right)dz_1.$$

Pour évaluer cette intégrale, on cherche à exprimer la forme
différentielle
$\left(\frac{1}{z_1-\bar{z}_1}-\frac{\partial}{\partial z_1}
h(z_1,\omega_c)\right)dz_1$ comme la différentielle totale d'une
fonction holomorphe. Commençons par noter que $h =
-4\Reel(\Lambda)$ est la partie réelle de la fonction
$-4\Lambda(z_1,\ldots,z_n)$ holomorphe par rapport à $z_1.$ On en
déduit que la forme $ \frac{\partial}{\partial z_1}
h(z_1,\omega_c)dz_1$ est la différentielle totale de la fonction
holomorphe $z_1\mapsto -2\Lambda(z_1,\omega_c).$

On veut faire de même avec la forme $\frac {dz_1}{z_1-\bar{z}_1}$
lorsque $z_1$ appartient au demi-cercle $\mc C$ de diamètre
$[\omega_{r_2},\omega_{r_1}].$ L'équation de $\mc C$ étant donnée
par
$$z_1\bar{z}_1- \frac{\omega_{r_1}+\omega_{r_2}} 2
(z_1+\bar{z}_1)+\omega_{r_1}\omega_{r_2}=0,$$ on trouve d'abord
\begin{equation*}\label{primlog}\frac 1{z_1-\bar{z}_1}=
\frac 12\left(\frac 1 {z_1-\omega_{r_1}}+\frac 1
{z_1-\omega_{r_2}}\right).\end{equation*} On conclut que $$\frac
{dz_1}{z_1-\bar{z}_1}=\frac 12\,\mathrm{ d }
\Big(\ln(z_1-\omega_{r_1})+\ln(z_1-\omega_{r_2})\Big).$$

Nous obtenons ainsi l'égalité valable pour tout $\tau$ sur le
demi-cercle $\mc C$
\begin{equation}\label{ELambda1}
L_A^{(n-1)}(0)=\frac{(n-1)!\, 2^{n}R_F}{4i\pi
}\Big[4\Lambda(z_1,\omega_c)+\ln(z_1-\omega_{r_1})+\ln(z_1-\omega_{r_2})\Big]_{\tau}^{\tau^*}
.\end{equation}

On rappelle que pour tout $\tau\in \mc H$ on a l'égalité
(\ref{LAMBDATRANSFPHI}), c'est à dire
$$\Big[4\Lambda(z_1,\omega_c)\Big]_{\tau}^{\tau^*}=\ln
\left[-(c_1\tau+d_1)^2)\right] +
\sum_{k=2}^n\ln|c_k\omega_k+d_k|^2 +i\pi \Phi(A,\omega_c).$$ On
sait que $\epsilon=c\omega+d$ est une unité relative de $K/F.$ Ses
images complexes sont donc sur le cercle unité. Il s'ensuit que
$\ln|c_k\omega_k+d_k|=0$ pour $2\leq k \leq n.$ En outre, on a par
définition
$\Psi(A)=2^nR_F\Phi(A,\omega_c)-2^{n-2}R_F\sign(c_1\tra(A_1)).$

D'après l'identité (\ref{ELambda1}) et les remarques précédentes,
il nous suffit pour obtenir la formule souhaitée de montrer que
pour tout $\tau \in \mc H,$ l'égalité
\begin{equation}\label{Dz0calc}
\Big[\ln\left(z_1-\omega_{r_1}\right)+\ln\left(z_1-\omega_{r_2}\right)\Big]_{\tau}^{\tau^*}
=-\ln\left(-\left(c_1\tau+d_1\right)^2\right)
-i\pi\sign(c_1\tra(A_1))
\end{equation}
est satisfaite. Notons d'abord que le complexe
$$\tau^*-\omega_{r_1}=\frac{\tau-\omega_{r_1}}{\epsilon_{r_1}(c_1\tau+d_1)}$$
est un point de $\mc H$ et $(\tau-\omega_{r_1})^{-1}$ un point du
demi-plan inférieur. On en déduit que
$$\big[\ln\left(z_1-\omega_{r_1}\right)\big]_{\tau}^{\tau^*}=
-\ln\left(\epsilon_{r_1}(c_1\tau+d_1)\right).$$ En menant le même
calcul avec $\omega_{r_2}$ et
$\epsilon_{r_2}=\epsilon_{r_1}^{-1},$ on en conclut que
\begin{equation}\label{aideln}\Big[\ln(z_1-\omega_{r_1})+\ln(z_1-\omega_{r_2})\Big]_{\tau}^{\tau^*}\!=\!
-\ln\left(\epsilon_{r_1}(c_1\tau+d_1)\right)-
\ln\left(\epsilon_{r_2}(c_1\tau+d_1)\right).\end{equation}
L'égalité précédente va nous permettre d'établir (\ref{Dz0calc}).
Pour cela, distinguons deux cas selon le signe de $c_1\tra(A_1).$
C'est aussi le signe de $c_1\epsilon_{r_1}$ puisque l'on a
$\epsilon+\epsilon^{-1}=\tra(A).$
\begin{itemize}
\item Cas 1 : $c_1\tra(A_1)<0.$ Les complexes
$i\epsilon_{r_1}(c_1\tau+d_1)$ et $i\epsilon_{r_2}(c_1\tau+d_1)$
sont donc dans le demi-plan $\Reel(\tau)>0.$ Par suite
\end{itemize} \noindent
\begin{align*}\ln\left(-(c_1\tau+d_1)^2\right)-i\pi=&\ln\left(i\epsilon_{r_1}(c_1\tau+d_1)\right)+
\ln\left(i\epsilon_{r_2}(c_1\tau+d_1)\right)+2\ln(-i)
\\=&\ln\left(\epsilon_{r_1}(c_1\tau+d_1)\right)+
\ln\left(\epsilon_{r_2}(c_1\tau+d_1)\right).\end{align*}
\begin{itemize}
\item Cas 2 : $c_1\tra(A_1)>0.$ Les complexes
$-i\epsilon_{r_1}(c_1\tau+d_1)$ et $-i\epsilon_{r_2}(c_1\tau+d_1)$
sont dans le demi-plan $\Reel(\tau)>0.$ Par suite
\end{itemize}\noindent
\begin{align*}\ln\left(-(c_1\tau+d_1)^2\right)+i\pi=\ln\left(\epsilon_{r_1}(c_1\tau+d_1)\right)+
\ln\left(\epsilon_{r_2}(c_1\tau+d_1)\right).\end{align*}

Ces deux cas ont rassemblés dans la formule
\begin{equation*}\label{dert}\ln\left(-(c_1\tau+d_1)^2\right)+i\pi\sign(c_1\tra(A_1))=\ln\left(\epsilon_{r_1}(c_1\tau+d_1)\right)+
\ln\left(\epsilon_{r_2}(c_1\tau+d_1)\right).\end{equation*} En
combinant l'égalité précédente et (\ref{aideln}), on obtient
l'identité (\ref{Dz0calc}) souhaitée. Le théorème s'ensuit.
\end{proof}









\subsection{Propriétés de l'invariant $\Psi(A)$. }\label{calc}
La définition \ref{defPsi} nous permet d'associer à une matrice
$A\in \Gamma$ lorsqu'elle est elliptique une famille de $n$
nombres réels $\{\Psi_j(A),1\leq j \leq n\},$ et lorsqu'elle est
quasi-elliptique un unique réel $\Psi_j=\Psi(A).$

La proposition \ref{propinvar} rassemble les propriétés communes
de ces invariants. Elle généralise à $\Gamma$ les résultats de
Rademacher [Ra1, Satz 7 et 9] pour $SL_2(\mb Z).$ Nous commençons
par montrer que dans le cas elliptique on peut calculer
$\Psi_j(A).$

\begin{lem}\label{lemmeellpsi}
Soit $A=\left(\begin{array}{cc} a &b\\
c& d
\end{array}\right)\in \Gamma$ une matrice elliptique. Alors $A$ est d'ordre
fini $m>1,$ et pour tout entier $j\in \{1,\ldots,n\},$ on a
$$\Psi_j(A)=-2^{n-2}R_F\left[\frac{\ln\left[-(c_j\omega_j+d_j)^2\right]}{i\pi}-\sign(c_j\tra(A_j))\right],$$ où $\omega_j$ est l'unique point fixe de
$A_j$ dans $\mc H,$ et $c_j\omega_j+d_j$ est une valeur propre de
$A_j.$

En particulier, $2\Psi_j(A)/R_F$ est un rationnel dont le
dénominateur divise $m.$

\end{lem}
\begin{proof} La preuve est immédiate.
Pour tout entier $k,$ $1\leq k\leq n,$ la composante $A_k$ de la
matrice $A$ a un unique point fixe $\omega_k$ dans $\mc H.$ On
déduit de (\ref{LAMBDATRANSFPHI}) que
$$\ln\left[-(c_j\omega_j+d_j)^2\right]+\sum_{k\neq
j}\ln|c_k\omega_k+d_k|^2+4i\pi\Phi_j(A,(\omega_k)_{k\neq j})=0.$$
Puisque $A$ est elliptique, elle est d'ordre fini $m>1.$ Ainsi les
valeurs propres $c_k\omega_k+d_k$ de $A_k$ sont des racines
$m$-ièmes de l'unité, et l'égalité précédente se réduit à
$$\ln\left[-(c_j\omega_j+d_j)^2\right]+4i\pi\Phi_j(A,(\omega_k)_{k\neq j})=0.$$
Il suffit d'utiliser la définition \ref{defPsi} pour conclure.
\end{proof}

\smallskip
\begin{prop}\label{propinvar}
Soit $A\in \Gamma$ une matrice elliptique ou quasi-elliptique de
$\Gamma.$ On a les égalités
\begin{enumerate}
\item $$\Psi_j(P^{-1}AP)=\Psi_j(A)\ \ \ \ \ \ \textrm{pour toute
matrice } P \textrm{ de } \Gamma.$$

\item $$\Psi_j(-A)=\Psi_j(A),\qquad \mathrm{ et }\qquad
\Psi_j(A^{-1})=-\Psi_j(A).$$

\item Soit $A\in \Gamma$ une matrice quasi-elliptique, et $k\in \mb Z$. On a l'égalité
$$\Psi(A^k)=k\Psi(A).$$

\end{enumerate}
\end{prop}

\noindent \textit{Démonstration.} On note $\omega_c\in \mc
H^{n-1}$ le point fixe de $(A_k)_{k\neq j}.$ \noindent Commençons
par démontrer la deuxième assertion. L'identité
$\Psi_j(-A)=\Psi_j(A)$ résulte immédiatement de la définition
\ref{defPsi}. De plus, on obtient grâce au théorème \ref{phicomp}
l'égalité
\begin{equation*}\label{PhiAA-1}\Phi_j(A,\omega_c)+\Phi_j(A^{-1},\omega_c)=0.\end{equation*}
On en déduit que $\Psi_j(A^{-1})=-\Psi_j(A).$

\smallskip
\noindent On se donne maintenant une matrice $P$
de $\Gamma$ et on note $P^{-1}AP=\left(\begin{array}{cc} a' &b'\\
c'& d'
\end{array}\right)$.

Pour démontrer les assertions $1$ et $3,$ nous distinguons deux
cas selon que la composante $A_j$ de $A$ est elliptique ou
hyperbolique.
\begin{itemize}
\item[i)] Dans le premier cas, la matrice $A$ est elliptique de
point fixe $\omega\in \mc H^n.$ Par suite, $P^{-1}AP$ est
elliptique de point fixe $P^{-1}\omega\in \mc H^n.$ Un calcul
direct montre l'égalité suivante entre les valeurs propres de $A$
et de $P^{-1}AP$ :
\begin{equation*}c'_j(P^{-1}_j\omega_j)+d'_j=c_j\omega_j+d_j. \label{ctdPAP}\end{equation*}
On en déduit que $c_j$ et $c'_j$ ont même signe, et par suite
\begin{equation*} \sign(c'_j\tra(P^{-1}_jA_jP_j))=\sign(c_j\tra(A_j)).\label{sgnecPAP}\end{equation*}

En utilisant les deux égalités précédentes ainsi que le lemme
\ref{lemmeellpsi}, on obtient le résultat souhaité :
$\Psi_j(P^{-1}AP)=\Psi_j(A).$

\item[ii)] Supposons maintenant $A_j=A_1$ hyperbolique. Le
théorème \ref{thL0PSI} relie $\Psi(A)$ à $ L_A^{(n-1)}(0).$ La
proposition \ref{propLAPAP} permet d'établir facilement les
propriétés demandées : \end{itemize}$$\Psi(P^{-1}AP)=\Psi(A)
\qquad \mathrm{et} \qquad \Psi(A^k)=k\Psi(A).$$ \hfill $\Box$

\medskip
Si la matrice $A$ de $ \Gamma$ est elliptique, les invariants
$\{\Psi_j(A)/R_F,1\leq j\leq n\}$ sont rationnels d'après le lemme
\ref{lemmeellpsi}. Si $A$ est quasi-elliptique, nous allons voir
que l'invariant $\Psi(A),$ bien que défini de manière analogue, a
une nature très différente dès que $n>1.$

On se propose dans la partie suivante de souligner l'intérêt
arithmétique des invariants $\Psi(A)$ si $A$ est quasi-elliptique.
On verra plus particulièrement que pour $n=2,$ la nature
arithmétique de ces invariants de classes de $\Gamma$ est
gouvernée par la conjecture de Stark. Dans les rares cas où cette
conjecture est démontrée, cela nous permettra de calculer quelques
invariants associés à des matrices quasi-elliptiques.


\subsection{Lien entre l'invariant $\Psi(A)$ et la conjecture de
Stark.}\label{PsiStark} Soit $F$ un corps quadratique réel de
nombre de classes 1, et $K$ une extension quadratique
quasi-totalement complexe. On note $v_1,$ $v_2$ les deux places
réelles de $K$ au-dessus de $\iota_1$ et $v_c$ la place complexe
de $K$ au-dessus de $\iota_2.$ On note $H$ (resp. $H^+$) le corps
de classes de Hilbert (resp. au sens restreint) de $K,$ et $G$ le
groupe de Galois de l'extension $H^+/K.$

On suppose que $K$ n'a pas d'unité de norme $-1.$ Par conséquent,
l'extension $H^+/H$ est une extension quadratique, et on note
$\rho\in G$ le générateur de son groupe de Galois. Le groupe $G$
est donc d'ordre $2h_K,$ où $h_K$ désigne le nombre de classes de
$K.$

Nous allons associer à chaque élément $\sigma$ du groupe $G$ une
classe de conjugaison $\mf C_\sigma$ de $\Gamma$ constituée de
matrices quasi-elliptiques $A_\sigma.$

Le but de ce paragraphe est de formuler précisément la conjecture
de Stark associée à l'extension $H^+/K$ pour prédire la nature
arithmétique des valeurs spéciales $L_{A_\sigma}'(0)$ et par suite
de la famille d'invariants de classes $\{\Psi(\mf C_\sigma), \,
\sigma\in G \}.$

L'application d'Artin permet d'associer à tout idéal fractionnaire
$\mf a$ de $K$ un élément $\sigma_\mf a$ de $G.$ D'après la
théorie du corps de classes, cette application induit un
isomorphisme du groupe des classes de $K$ au sens restreint
$\Cl^+_K$ sur le groupe $G.$ Notons $S=\{v_{1},v_2,v_c\}$
l'ensemble des places archimédiennes de $K.$ On associe à tout
élément $\sigma$ de $G$ la fonction zeta partielle définie pour
$\Reel(s)>1$ par
$$\zeta_S(s,\sigma):=\sum_{\mf m \subset \mc O_K \atop \sigma_{\mf
m}=\sigma} N(\mf m)^{-s},$$ la sommation portant sur les idéaux
entiers $\mf m$ de $K$ dont le symbole d'Artin $\sigma_\mf m$ est
égal à $\sigma.$ Ces fonctions admettent un prolongement
méromorphe pour $s\in \mb C.$

Avant de pouvoir énoncer la conjecture de Stark, nous devons
étudier l'ordre d'annulation de la fonction $\zeta_S(s,\sigma)$ en
$s=0.$ Les résultats souhaités se déduisent du comportement en
$s=0$ des fonctions $L$ de Hecke de l'extension $H^+/K$ définies
comme suit. On associe à tout caractère $\chi$ sur le groupe
abélien $G$ la fonction $L$ de Hecke définie pour $\Reel(s)>1$ par
$$L_S(s,\chi)=\sum_{\sigma\in G}\chi(\sigma)\zeta_S(s,\sigma).$$
Ces fonctions ont un développement en produit d'Euler. Si $\chi$
n'est pas le caractère trivial, elle admettent un prolongement
holomorphe à $\mb C.$

L'ordre d'annulation $r(\chi)$ de $L_S(s,\chi)$ en $s=0$ est
connu. Il se déduit de la propriété $L_S(1,\chi)\neq 0$ et de
l'équation fonctionnelle de $L_S(s,\chi)$ (cf. [We]). On obtient
les résultats généraux suivants ([D-S-T]) : si $\chi$ est le
caractère trivial, on a $r(\chi)=|S|-1;$ sinon, $r(\chi)$ est le
nombre de places $v$ de $S$ dont le groupe de décomposition
$D_v\subset G$ est contenu dans le noyau de $\chi.$ Sous nos
hypothèses, l'extension $H^+/H$ est quadratique. Par conséquent,
la seule place de $S$ qui se décompose totalement dans l'extension
$H^+/K$ est $v_c,$ tandis que les deux autres places $v_1$ et
$v_2$ se ramifient. On en déduit que
\begin{itemize}

\item $r(\chi)=2$ si $\chi$ est le caractère trivial,

\item $r(\chi)=1$ si $\chi(\varrho)=-1,$

\item $r(\chi)=3$ sinon.
\end{itemize}
Les formules d'orthogonalité des caractères permettent de relier
$\zeta_S(s,\sigma)$ aux fonctions $L_S(s,\chi)$ selon la formule
\begin{equation*}\zeta_S(s,\sigma)=\frac{1}{2h_K}\sum_{\chi\in
\hat{G}}\overline{\chi}(\sigma)L_S(s,\chi).\end{equation*} On en
conclut finalement que les fonctions $\zeta_S(s,\sigma)$
s'annulent en $s=0$ pour tout $\sigma$ de $G.$

On obtient même un renseignement supplémentaire : l'égalité
précédente permet aussi d'établir l'identité
$$\zeta_S(s,\sigma)+\zeta_S(s,\sigma\varrho)=\frac{1}{h_K}\sum_{\chi\in
\hat{G}\atop \chi(\varrho)=1}\overline{\chi}(\sigma)L_S(s,\chi),$$
la sommation portant sur les caractères $\chi$ de $G$ tels que
$\chi(\varrho)=1.$ Toutes les fonctions $L$ qui apparaissent dans
cette dernière somme vérifient $L'_S(0,\chi)=0$ d'après la
discussion précédente. On en déduit que pour tout $\sigma \in G$
on a l'égalité
$$\zeta'_S(0,\sigma)+\zeta_S'(0,\sigma\varrho)=0.$$

La conjecture de Stark s'énonce de la manière suivante dans notre
contexte ([Ta],[D-S-T]).
\medskip

\begin{conj}[Stark]~
Soit $q$ le nombre de racines de l'unité contenues dans $H^+.$ Il
existe une unité de $u$ de $H^+$ qui vérifie les propriétés qui
suivent.
\begin{itemize}
\item[a)] Soit $w'$ une place de $H^+$ qui ne divise pas
l'unique place complexe $v_c$ de $K.$ Alors on a l'égalité
$|u|_{w'}=1.$

\item[b)] Soit $w$ une place de $H^+$ au-dessus de $v_c.$ Alors pour tout $\sigma\in
G$ on a l'égalité $$\ln |u^\sigma|_w=-q\,\zeta_S'(0,\sigma).$$
\item[c)] L'extension $K(u^{1/q})$ est une extension abélienne de $K.$
\end{itemize}

\end{conj}

Nous reformulons maintenant cette conjecture en termes de valeurs
spéciales $L_A'(0),$ où $A$ est une matrice quasi-elliptique de
$\Gamma.$ Pour cela, nous associons à tout élément $\sigma$ de $G$
une classe de conjugaison $\mf C_\sigma$ de $\Gamma$ selon la
construction suivante.

On commence par associer à $\sigma\in G$ un idéal fractionnaire
$\mf b$ de $K$ tel que $\sigma_\mf b=(\sigma)^{-1}.$ Considéré
comme un $\mc O_F$-module libre de rang 2, l'idéal $\mf b$ a une
base $(\gamma,$ $\widetilde{\gamma})\in K^2.$ On pose
$\omega:=\widetilde{\gamma}/\gamma.$ Par convention, la base
$(\gamma,\widetilde{\gamma})$ de $\mf b$ est orientée : on demande
que $\omega$ vérifie $\omega_{r_2}<\omega_{r_1},$ où
$\alpha_{r_1}$ et $\alpha_{r_2}$ désignent les images réelles d'un
élément de $\alpha\in K$ associées aux places $v_1$ et $v_2.$
L'idéal fractionnaire $\mf a:=(\gamma^{-1})\mf b$ s'écrit $\mf
a=\mc O_F +\omega\mc O_F$ et vérifie $\sigma_\mf
b=\sigma_{(\gamma)}\sigma_\mf a.$

Ceci étant, on remarque que le groupe des unités relatives de
l'extension $K/F$ est de rang 1. On note $\epsilon$ le générateur
de sa partie libre vérifiant $\epsilon_{r_1}>1.$ On a $\epsilon\mf
a=\mf a.$ Il existe donc des entiers $a,$ $b,$ $c,$ $d$ de $\mc
O_F$ tels que $\epsilon\omega=a\omega+b$ et $\epsilon=c\omega+d.$
Autrement dit, on a l'égalité matricielle
\begin{equation}\label{relmatomega}\left(\begin{array}{cc} \omega &\omega^g\\
1& 1
\end{array}\right)^{-1}\left(\begin{array}{cc} a &b\\ c& d
\end{array}\right)\left(\begin{array}{cc} \omega &\omega^g\\ 1&1
\end{array}\right)=\left(\begin{array}{cc} \epsilon &0\\ 0& \epsilon^g
\end{array}\right),\end{equation}
où on a noté $g$ le générateur du groupe de Galois de l'extension
$K/F.$
On en déduit que $ad-bc=\epsilon\epsilon^g=1.$ Par suite, la matrice $A=\left(\begin{array}{cc} a &b\\
c& d
\end{array}\right)$ est un élément de $\Gamma$ quasi-elliptique de points
fixes $\{\omega,\omega^g\}.$

La matrice $A$ construite précédemment dépend de l'élément
$\sigma$ de $G$ mais aussi du choix arbitraire de l'idéal $\mf b$
tel que $\sigma_\mf b=(\sigma)^{-1}$ et d'une $\mc O_F$-base
orientée de $\mf b.$ On montre facilement que des choix différents
conduisent à une matrice conjuguée de $A$ dans $\Gamma.$ Par
suite, $\sigma$ définit une classe de conjugaison $\mf C_\sigma$
de $\Gamma.$

Considérons la fonction $L_A(s)$ associée à la matrice $A.$ Avec
les notations du paragraphe \ref{partiequasiell}, le $\mc
O_F$-module $\mc M_\omega=\mc O_F+\omega\mc O_F$ n'est autre que
l'idéal fractionnaire $\mf a.$ Le groupe des unités $U_A$ est un
sous-groupe d'indice $e_K=2$ ou $1$ du groupe des unités de $K,$
selon que l'unité fondamentale de $F$ est ou n'est pas la norme
relative d'une unité de $K.$ En outre, l'hypothèse
$\epsilon_{r_1}>1$ entraîne que $\tra(A_1)>2.$ En utilisant
(\ref{signcepsi}), on en déduit que $c_1\tra(A_1)>0.$ Par
conséquent, la fonction $L_A$ associée à $A$ est donnée d'après la
définition \ref{defLA} par la formule
$$L_A(s)=e_K\sum_{(\alpha)\subset \mf a}\frac{\varphi(\alpha)}{|N_{K/\mb
Q}(\alpha)|^s},$$ où on a noté $\varphi$ le caractère de $K^*$
défini par $\varphi(\alpha)=\sign(\alpha_{r_1}\alpha_{r_2}).$

Nous relions maintenant $L_A(s)$ à la fonction
$\zeta_S(s,\sigma).$ Soit $\mf m$ un idéal de $\mc O_K$ tel que
$\sigma_\mf m=\sigma.$ On a alors $\sigma_{\mf b \mf m}=1.$
L'idéal $\mf b\mf m$ est donc un idéal principal $(\beta )$
engendré par un élément $\beta\in \mf b$ totalement positif. On en
déduit que
\begin{align*}\zeta_S(s,\sigma)=N(\mf b
)^{s}\sum_{\mf m\subset \mc O_K\atop \sigma_\mf m=\sigma}N(\mf b
\mf m)^{-s} \\= N(\mf b )^{s}\sum_{(\beta)\subset \mf b \atop
\beta \gg 0} |N_{K/\mb Q}(\beta)|^{-s}.
\end{align*}
Soit $\mf m'$ un idéal de $\mc O_K$ tel que $\sigma_{\mf
m'}=\sigma\varrho.$ L'idéal $\mf m'\mf b$ est un idéal principal
$(\beta)$ engendré par un élément $\beta\in \mf b$ tel que
$\varphi(\beta)=-1.$ Il s'ensuit que
$$\zeta_S(s,\sigma)-\zeta_S(s,\sigma\varrho)=N(\mf
b)^s\sum_{(\beta)\subset \mf b}\frac{\varphi(\beta)}{|N_{K/\mb
Q}(\beta)|^s}.$$ On effectue le changement de variable
$\alpha:=\gamma^{-1}\beta$ pour obtenir
\begin{equation}\label{lienLAzeta}L_A(s)=e_K\varphi(\gamma)N(\mf
a)^{-s}\left[\zeta_S(s,\sigma)-\zeta_S(s,\sigma\varrho)\right].\end{equation}

On conclut finalement de l'égalité
$\zeta_S'(0,\sigma\varrho)=-\zeta_S'(0,\sigma)$ que

\smallskip
\begin{equation}\label{LZETAPART}L'_A(0)=\pm 2e_K\zeta_S'(0,\sigma),\end{equation}

\smallskip
\noindent le signe étant $+1$ ou $-1$ selon que $\sigma_{\mc
M_\omega}^{-1}$ est égal à $\sigma$ ou $\sigma\varrho$
respectivement. On remarque que le réel $$L'_A(0)=\Psi(A)$$ ne
dépend que de la classe de conjugaison $\mf C_\sigma$ de $A$
d'après la proposition \ref{propinvar}.3.

\smallskip
Soit $(A_\sigma)_{\sigma\in G}\in \Gamma^{2h_K}$ une famille de
représentants des classes de conjugaison $(\mf
C_\sigma)_{\sigma\in G}$ construites précédemment. La conjecture
de Stark prédit donc que le réel $e^{-\frac
{q}{2e_K}\Psi(A_\sigma)}$ est la valeur absolue d'une unité
algébrique qui engendre l'extension $H^+/K.$ Il faut remarquer que
l'on a reformulé une conjecture portant sur la valeur spéciale de
fonctions $L$ en $s=0$ en une conjecture portant sur la valeur
spéciale de la fonction réelle-analytique $\Phi(A_\sigma,z_2)$ en
un point algébrique $\omega_c\in\mc H$ associé à $A_\sigma.$ La
conjecture de Stark est démontrée dans le cas ``trivial'' où
$[H^+:K]=2$ (voir [Ta, Th. 5.4]). On dispose alors d'une formule
explicite pour $L'_{A_\sigma}(0)$ que nous utilisons dans le
paragraphe suivant pour obtenir la valeur de certains invariants
$\Psi(A).$

\subsection{Exemples numériques.}
Soit $F$ un corps de nombres totalement réel de degré $n,$ de
nombre de classes 1. Soit $K$ une extension quadratique de $F$
quasi-totalement complexe. On note $H^+$ le corps de classes de
Hilbert au sent restreint de $K.$

Le but de ce paragraphe est de calculer explicitement un certain
invariant $\Psi(A)$ sous l'hypothèse que l'extension $H^+/K$ est
quadratique. On profite en effet de cette hypothèse très
restrictive pour calculer $L'_A(0)$ et donc $\Psi(A).$

Soit $\chi$ le caractère non-trivial du groupe des classes au sens
restreint de $K.$ Dans cette situation, on a l'égalité
\begin{equation}\label{Lquotzeta}L(s,\chi)=\frac{\zeta_{H^+}(s)}{\zeta_K(s)}.\end{equation}
Nous pouvons calculer le premier terme non nul du développement en
série de Laurent au voisinage de $s=0$ de ce quotient.

Supposons pour simplifier que l'extension $F(i)/F$ est non
ramifiée en 2. Il s'ensuit que le corps $K(i)$ est contenu dans
$H^+.$ Par hypothèse, $H^+$ est une extension quadratique de $K;$
on en déduit l'égalité $$H^+=K(i).$$

Il y a trois extensions quadratiques de $F$ contenues dans $K(i)$
: le corps $F(i)$ qui est CM, le corps quasi-totalement complexe
$K,$ et un corps quasi-totalement réel noté $K'.$ La fonctorialité
des fonctions $L$ d'Artin permet d'écrire que
\begin{equation}\label{quotzeta}\frac{\zeta_{K(i)}}{\zeta_F}=\frac{\zeta_{F(i)}}{\zeta_F}
\frac{\zeta_{K}}{\zeta_F}\frac{\zeta_{K'}}{\zeta_F}.\end{equation}
On rappelle que la fonction zeta d'un corps de nombres $k$ a le
développement de Taylor suivant au voisinage de $s=0$ :
$$\zeta_k(s)=-\frac{h_k R_k}{W_k}s^{u+v-1}+O(s^{u+v}),$$ où $W_k$ est le nombre
de racines de l'unités de $k,$ $u$ (resp. $v$) le nombre de places
réelles (resp. complexes) de $k$ , $h_k$ le nombre de classes de
$k$ et $R_k$ le régulateur de $k.$ On déduit de (\ref{Lquotzeta}),
(\ref{quotzeta}) et de l'égalité précédente que
\begin{equation}\label{LLautres}
L^{(n-1)}(0,\chi)=(n-1)!\frac{2h_{F(i)}h_{K'}}{W_{F(i)}}\frac{R_{K'}}{R_F}.\end{equation}

Notons que si $n=2$, cette identité montre que la conjecture de
Stark est valide dans la situation ``triviale'' où $H^+=K(i).$ En
effet, on établit facilement que le réel $R_{K'}/R_F$ est le
logarithme d'une unité de $K'.$

\medskip

Revenons au calcul des invariants de classes de $\Gamma.$ La
fonction $L(s,\chi)$ peut s'exprimer en fonction de $L_A(s)$ pour
une certaine matrice quasi-elliptique $A$ de $\Gamma.$ Pour ce
faire, on considère l'anneau des entiers $\mc O_K$ de $K$ comme un
$\mc O_F$ module libre de rang 2. Il a une base orientée
$(\gamma,\tilde{\gamma})\in K^2.$ Cette base induit alors une
matrice quasi-elliptique $A$ de $\Gamma$ en suivant une
construction similaire à celle du paragraphe \ref{PsiStark}. On
déduit de (\ref{lienLAzeta}) l'égalité
\begin{equation*}\label{LALCHI}L_A(s)=e_K N_{K/\mb Q}(\gamma)^s
\varphi(\gamma)L(s,\chi).\end{equation*} A l'aide de
(\ref{LLautres}) et du théorème \ref{thL0PSI}, on en conclut
finalement que

\begin{equation}\label{PSIACALC}
\Psi(A)=\varphi(\gamma)\frac{2e_Kh_{F(i)}h_{K'}R_{K'}}{W_{F(i)}R_F}.
\end{equation}

\noindent \textbf{Exemples numériques.} On se limite au cas où
$n=2.$ Le corps $F=\mb Q(\sqrt{7})$ est de nombre de classes 1.
Son unité fondamentale est $u=8+3\sqrt{7}.$ Le corps $F(i)$ est de
nombre de classes 1, et il contient 4 racines de l'unité.

On considère d'abord l'extension quadratique $K$ de $F$ engendrée
par le réel $\omega_{r_1}^K=\sqrt{-2+\sqrt{7}}.$ L'anneau des
entiers de $K$ a pour base relative orientée $(1,\,\omega^K).$ Le
groupe des unités de $K$ est de rang 2, engendré par
$\{-1,u,\epsilon^K\}$ où $\epsilon^K$ est de norme relative 1 et a
pour image réelle
$\epsilon_{r_1}^K:=(3+\sqrt{7})\omega_{r_1}^K-2-\sqrt{7}.$ On
construit alors grâce à la relation (\ref{relmatomega}) la matrice
$$A_1:=\left(\begin{array}{cc} -2-\sqrt{7} &1+\sqrt{7}\\
3+\sqrt{7}& -2-\sqrt{7}
\end{array}\right),$$ image par le plongement $\iota_1$ d'une matrice quasi-elliptique $A$ de $SL_2(\mc
O_F).$ On a $h_K=e_K=1.$ Le corps de classes de Hilbert au sens
restreint de $K$ est $K(i).$ Les trois extensions de $F$ contenues
dans $K(i)$ sont $K,$ $F(i)$ et le corps $K'$ engendré par le réel
$\omega_{r_1}^{K'}=\sqrt{2+\sqrt{7}}.$ Le groupe des unités de
$K'$ a pour générateurs $\{-1,u,\epsilon^{K'}\}$ où
$\epsilon_{r_1}^{K'}:=(9+3\sqrt{7})\omega_{r_1}^{K'}+18+7\sqrt{7}$
et $N_{K'/F}(\epsilon^{K'})=1.$ Par suite, un calcul élémentaire
montre que $R_{K'}/R_F=2\ln \epsilon_{r_1}^{K'}.$ En outre, on
$h_{K'}=e_{K'}=1.$ On déduit donc de l'égalité (\ref{PSIACALC})
que l'invariant associé à la matrice $A$ est

$$\Psi(A)=\ln \epsilon_{r_1}^{K'}=\ln\left(\left(9+3\sqrt{7}\right)\sqrt{2+\sqrt{7}}+18+7\sqrt{7}\right).$$

Nous pouvons exprimer l'égalité ci-dessus en termes de valeur
spéciale en $\omega_c^K=i\sqrt{2+\sqrt{7}}\in \mc H$ de la somme
de Dedekind généralisée $s(d,c;\,.),$ où $d_1=-2-\sqrt{7},$
$c_1=3+\sqrt{7}.$ Grâce à la proposition \ref{remeucl}, on peut
même se ramener à la somme de Dedekind fondamentale
$s=s(0,1;\,.).$ A l'aide de l'égalité (\ref{sdcs01exple}), on
obtient ainsi que $\Psi(A)/R_F$ est égal à
$$4s\left(\frac{\omega_c^K} {-\omega_c^K+1}\right)-4s\left(\omega_c^K+1\right)+2+\kappa_F\sqrt{2+\sqrt{7}}
\left[\frac{72037+23827\sqrt{7}}{4683}\right].$$

On remplace maintenant le corps $K$ par $K'.$ La construction
précédente conduit à une matrice hyperbolique
$$A_1'=\left(\begin{array}{cc} 18+7\sqrt{7} &39+15\sqrt{7}\\
9+3\sqrt{7}&18 +7\sqrt{7}
\end{array}\right)$$ dont les points fixes réels sont $\pm \sqrt{2+\sqrt{7}}=\pm\omega_{r_1}^{K'}.$ En échangeant les rôles de $K$ et $K',$
nous pouvons calculer l'invariant $\Psi(A').$ On trouve

\begin{equation*}\Psi(A')=\ln\epsilon_{r_1}^K=\ln
\left(\left(3+\sqrt{7}\right)\sqrt{-2+\sqrt{7}}-2-\sqrt{7}\right).\end{equation*}

\section{Analogies entre la conjecture de Stark et celle de Darmon.}
On se limite \`a $n=1$ ou 2. On se propose de souligner les
analogies entre la construction des points de Heegner  propos\'ee
par Darmon dans [Da] et notre construction bas\'ee sur les
s\'eries d'Eisenstein de poids 2 pour $\Gamma$.

Pour  $N\in \mb N$ sans facteur carr\'e,  on note $\Gamma_0(N)$ l'ensemble  des matrices $\left(\begin{array}{cc} a &b\\
c& d
\end{array}\right)$ de $\Gamma$ telles que $N|c.$ On se donne une forme
modulaire $f$ propre pour les op\'erateurs de Hecke, normalis\'ee,
de poids parall\`ele 2 sur $\Gamma_0(N).$

\subsection{Cas o\`u $F=\mb Q.$} Rappelons tout d'abord tr\`es
sommairement la construction classique des points de Heegner d'une
extension quadratique imaginaire de $\mb Q$ ({Da, Part. 2 et 3]).
 Si $f$ est cuspidale, on peut lui associer  une courbe
elliptique $\mathcal{E}$ d\'efinie sur $\mb Q.$  Notons
$\Lambda_{\mathcal{E}}\subset \mb C$ le r\'eseau de N\'eron de
$\mathcal{E}$, $\wp$ la fonction de Weierstrass associ\'ee. Pour
$\tau\in\mc H,$ on note $\Phi_N(\tau)\in \mathcal{E}(\mb C)$
l'image de
$$z_\tau:=2i\pi\int_{i\infty}^ \tau f(z)dz$$ par l'application
$(\wp,\wp') : \mb C/\Lambda_{\mathcal{E}} \rightarrow
\mathcal{E}(\mb C).$

 Si $\tau\in \mc H$ est un point
quadratique imaginaire, alors $\Phi_N(\tau)$ est alg\'ebrique:
c'est un point de $\mathcal{E}(H),$ o\`u $H$ est le corps de
classes de $\mb Q(\tau)$ associ\'e \`a un certain ordre $\mc
O_\tau$ ([Da, Th. 3.6]).

\smallskip

 Que se passe-t-il quand on remplace la forme cuspidale
$f$ de poids 2 par la s\'erie d'Eisenstein $E_{\mb Q,2}(z)$ de
poids 2? Cette s\'erie est la valeur en $s=0$ de la fonction
d\'efinie pour $\Reel(s)>1$ par $$ E_{\mb
Q,2}(z,s):=\sum_{(m,n)\in \mb Z^ 2}\!\!\!\prim (mz+n)^ {-2}y^
s|mz+n|^ {2s}=\frac {4i} s\frac
\partial{\partial z}E_\mb Q(z,s+1).$$ Les r\'esultats suivants se
d\'eduisent facilement de la formule limite de Kronecker (Th.
\ref{propEasai}.5) : $\frac 1{4\pi^2}\left( E_{\mb Q,2}(z)+\frac
\pi y\right)=-1/24+\ldots$ est une fonction holomorphe dont les
coefficients de Fourier sont rationnels; On a l'\'equation
fonctionnelle $E_{\mb Q,2}(Az)=(cz+d)^2E_{\mb Q,2}(z),$  et \`a
une constante d'int\'egration pr\`es on obtient finalement
l'identit\'e
\begin{equation}\label{Eetaint}2 \ln
\eta(\tau)= 2i\pi\int^\tau \frac 1 {4\pi^2}\left(E_{\mb
Q,2}(z)+\frac \pi y\right)dz.\end{equation} Noter que le terme
correctif $\frac \pi y $ a \'et\'e introduit pour obtenir une
forme diff\'erentielle ferm\'ee.

 En outre, la fonction $\eta(\tau)$ prend
des valeurs alg\'ebriques quand $\tau$ est quadratique imaginaire;
l'\'egalit\'e (\ref{Eetaint}) est donc un ingr\'edient crucial de
la d\'emonstration de  la conjecture de Stark pour le corps $\mb
Q(\tau).$

La morale est donc la suivante : on se donne  $f$ une forme
modulaire propre normalis\'ee sur de poids 2 pour $\Gamma_0(N)$,
et $\omega_f:=2i\pi f(z)dz$ la forme diff\'erentielle
$\Gamma_0(N)$-invariante associ\'ee. Soit $\tau \in\mc H$ un point
quadratique imaginaire. Alors $\int^\tau \omega_f$ permet de
construire des points alg\'ebriques sur la courbe elliptique $\mc
E$ si $f$ est une forme cuspidale, et des logarithmes de nombres
alg\'ebriques reli\'es \`a la conjecture de Stark si $f$ est la
s\'erie d'Eisenstein.

\medskip

\subsection{Cas o\`u $F$ est quadratique r\'eel de nombre de
classes au sens restreint $h_F^+=1$.}

On r\'esume sommairement la construction de Darmon. On se donne
$\mc E$ une courbe elliptique sur $F$ de conducteur $N$ associ\'ee
\`a une forme modulaire de Hilbert cuspidale $f$ de poids (2,2).

 Soit $\epsilon\in F$ une unit\'e telle
que $\epsilon_1>0$ et $\epsilon_2<0.$ On d\'efinit la forme
diff\'erentielle
$$\omega_f^+:=-\frac{4\pi^2}
{d_F}\left[f(z_1,z_2)dz_1dz_2+f(\epsilon_1z_1,\epsilon_2\bar{z}_2)d(\epsilon_1z_1)d(\epsilon_2\bar{z}_2)\right].$$
Soit $A\in \Gamma$ une matrice quasi-elliptique, et $\tau_A\in \mc
H$ le point fixe de sa composante elliptique $A_1$. Soit $K$
l'extension quasi-elliptique associ\'ee.

On choisit arbitrairement un point $x\in \mc H$ et on pose
$$J_{\tau_A}=\int ^{\tau_A}\int_x^{A_2x} \omega_f^+.$$ Ceci ne
d\'epend pas du choix de $x.$  Il existe conjecturalement un
r\'eseau $\Lambda_\mc E$ qui ne d\'epend que de $\mc E$ tel que,
en notant $\pi : \mb C/\Lambda_\mc E \rightarrow \mc E(\mb C)$
l'uniformisation de Weierstrass et  $q_{\mc E(K)}$ le cardinal du
groupe de torsion de $\mc E(K),$ le point $q_{\mc E(K)}
\pi\left(J_{\tau_A}\right)$ devrait \^etre un point de $\mc E(H),$
o\`u $H$ est un corps de classes pr\'ecis de $K.$

\smallskip

 Remplaçons maintenant la forme cuspidale $f$ de la
construction de Darmon par la s\'erie d'Eisenstein
$E_{F,2}(z_1,z_2)$ de poids 2 de $\Gamma.$ Elle est donn\'ee par
la valeur en $s=0$ de la fonction
 d\'efinie pour $\Reel(s)>1$ par
\begin{align*}E_{F,2}(z_1,z_2,s)=&\!\!\!\!\sum_{(\mu,\nu)\in \mc
O_F^2/U_F}\!\!\!\!\!\!\prim
 \frac 1{(\mu_1z_1+\nu_1)^{2}(\mu_2z_2+\nu_2)^{2}}
 \frac{ \left(y_1y_2\right)^s }{|\mu_1z_1+\nu_1|^{2s}|\mu_2z_2+\nu_2|^{2s}}  \\=& -\frac {4} {s^2}
 \frac \partial{\partial z_1}\frac \partial{\partial
 z_2}E_F(z_1,z_2,s+1).\end{align*}

On d\'eduit imm\'ediatement de  la formule limite de Kronecker
g\'en\'eralis\'ee que $\pi^2 E_{F,2}(z_1,z_2)$ est une vraie forme
modulaire de Hilbert holomorphe de poids (2,2) pour $\Gamma.$ Ses
coefficients de Fourier sont rationnels,  et on a l'\'egalit\'e
\begin{equation*}
E_{F,2}(z_1,z_2)=-\frac {4\pi^2}{d_F}\frac
\partial {\partial z_1}\frac \partial {\partial z_2} h(z_1,z_2).
\end{equation*}

Par suite, la forme diff\'erentielle associ\'ee \`a $E_{F,2}$ est
$$\omega_{E_{F,2}}^+=\frac
\partial {\partial z_1}\frac \partial {\partial z_2} h(z_1,z_2)dz_1dz_2+\frac
\partial {\partial z_1}\frac \partial {\partial z_2} (h)
(\epsilon
z_1,\epsilon_2\bar{z}_2)d(\epsilon_1z_1)d(\epsilon_2\bar{z}_2).$$

Supposons cette fois que  $A_1$ est hyperbolique et $A_2$
elliptique. Il r\'esulte alors de la formule (\ref{ELambda1}) que,
\`a une constante d'int\'egration pr\`es, on a  l'identit\'e
$$\frac {2R_F}{i\pi}\left[\int_x^{A_1x}\int^{\tau_A}
\omega_{E_{F,2}}^+ - \int_x^{A_1x} \frac 12\left[\frac 1
{z_1-\omega_{r_1}}+\frac 1 {z-\omega_{r_2}}\right]dz_1\right]
=-L'_A(0).$$

 Notons  $q_K$ le nombre de
racines de l'unit\'e de $K.$ Si l'on en croit la conjecture de
Stark (cf. partie \ref{PsiStark}), $-q_K L'_A(0)$ devrait \^etre
le logarithme d'une unit\'e d'un corps de classes pr\'ecis $H$ de
$K.$

\noindent {\sc Adresse de l'auteur:}
\medskip

\noindent Institut de Math\'ematiques, Univ. Bordeaux~1, France
\smallskip

\noindent e-mail : pierre.charollois@math.u-bordeaux.fr

\end{document}